\documentclass[11 pt]{amsart}
\usepackage{amscd, amssymb}
\newtheorem{pr}{Proposition}
\newtheorem{lm}{Lemma}
\newtheorem{tm}{Theorem}

\newcommand{\proj}{\mathbf P}

\newcommand{\rarr}{\rightarrow}

\newcommand{\com}{\mathbb{C}}
\newcommand{\Q}{\mathbb{Q}}
\newcommand{\Z}{\mathbb{Z}}

\newcommand{\GG}{{\mathbf{g}}}
\newcommand{\NN}{{\mathbf{n}}}
\newcommand{\DD}{{\mathbf{d}}}
\newcommand{\mumu}{{\mathbf{\mu}}}
\newcommand{\M}{{\overline M}}
\newcommand{\len}{\ell}

\def\scup{\mathbin{\text{\scriptsize$\cup$}}}

\newcommand{\bpf}{\noindent {\em Proof.} }
\newcommand{\epf}{\qed \vspace{+10pt}}

\newcommand{\mata}{\mathbf{A}}
\newcommand{\matb}{\mathbf{B}}

\newcommand{\matc}{\mathbf{C}}

\newcommand{\im}{\text{Im}}
\newcommand{\Aut}{\text{Aut}}
\newcommand{\qbar}{\hat{q}}
\newcommand{\pbar}{\hat{p}}
\newcommand{\rbar}{\hat{r}}
\newcommand{\qti}{\tilde{q}}
\newcommand{\pti}{\tilde{p}}
\newcommand{\rti}{\tilde{r}}

\begin{document}
\title{Relative maps and tautological classes}
\author{C. Faber and R. Pandharipande}
\date{April 2003}
\maketitle

\pagestyle{plain}
\setcounter{section}{-1}
\section{Introduction}
\subsection{Tautological classes}

Let ${\overline{M}_{g,n}}$ be the moduli space of stable curves of
genus $g$ with $n$ marked points defined over ${\mathbb C}$. 
Let $A^*(\overline{M}_{g,n})$ denote the Chow ring (always
taken here with ${\mathbb Q}$-coefficients).
The system of
tautological rings is defined
to be the set of smallest $\Q$-subalgebras of the Chow rings,
$$R^*(\overline{M}_{g,n}) \subset A^*(\overline{M}_{g,n}),$$
satisfying the
following two properties:
\begin{enumerate}
\item[(i)] The system is closed under push-forward via
all maps forgetting markings:
$$\pi_*: R^*(\overline{M}_{g,n}) \rarr R^*(\overline{M}_{g,n-1}).$$
\item[(ii)] The system is closed under push-forward via
all gluing maps:
$$\iota_*: R^*(\overline{M}_{g_1,n_1\scup\{*\}})
\otimes_{\Q}
R^*(\overline{M}_{g_2,n_2\scup\{\bullet\}}) \rarr
R^*(\overline{M}_{g_1+g_2, n_1+n_2}),$$
$$\iota_*: R^*(\overline{M}_{g, n\scup\{*,\bullet\}}) \rarr
R^*(\overline{M}_{g+1, n}),$$
with attachments along the markings $*$ and $\bullet$.
\end{enumerate}
While the definition appears restrictive,
natural algebraic constructions typically yield Chow classes
lying in the tautological ring.  For example, the standard $\psi$,
$\kappa$, and $\lambda$ classes in $A^*(\overline{M}_{g,n})$ 
all lie in the tautological ring.
The tautological rings also possess a rich conjectural structure,
see [FP] for a detailed discussion.

The cotangent line classes $\psi$ are tautological by the
following construction.
For each marking $i$, let
$L_i$ denote the associated cotangent line bundle over
${\overline M}_{g,n}$. The class $\psi_i$ is the first Chern class
of $L_i$,
$$\psi_i = c_1(L_i) \in A^1(\overline{M}_{g,n}).$$
Let $\pi$ denote  the map forgetting the last marking,
$$\pi: \overline{M}_{g,n+1} \rarr \overline{M}_{g,n},$$
and
let $\iota$ denote the gluing map,
$$\iota: \overline{M}_{g,\{1,2,\ldots,i-1,*,i+1,\ldots, n\}} \times
\overline{M}_{0,\{\bullet,i,n+1\}} \longrightarrow \overline{M}_{g,n+1}.$$
The ${\mathbb Q}$-multiples of
the fundamental classes $[\overline{M}_{g,n}]$  
are
contained in the tautological rings (as ${\mathbb Q}$-multiples of the 
units in the subalgebras).
A direct calculation shows:
$$-\pi_*\Big( \big(\iota_*([\overline{M}_{g,n}] \times
[\overline{M}_{0,3}])\big)^2\Big) = \psi_i.$$
Hence, the cotangent line classes lie in the tautological rings.
A discussion of the  $\kappa$ and $\lambda$ classes 
can be found in [FP].

Gromov-Witten theory defines natural classes in $A^*(\overline{M}_{g,n})$.
Let $X$ be a nonsingular projective variety, and let
$\overline{M}_{g,n}(X,\beta)$ be the moduli space of stable
maps representing $\beta\in H_2(X,{\mathbb Z})$.
Let $\rho$ denote the map to the moduli of curves,
$$\rho: \overline{M}_{g,n}(X,\beta) \rarr {\overline M}_{g,n},$$
in case $2g-2+n>0$.
Let $\omega \in A^*(\overline{M}_{g,n}(X,\beta))$ be
a Gromov-Witten class composed of algebraic primary fields and 
descendents. Then,
$$\rho_*\Big( \omega \cap [\overline{M}_{g,n}(X,\beta)]^{vir}\Big) \in 
A^*(\overline{M}_{g,n}).$$
The push-forwards of all Gromov-Witten classes of compact
homogeneous varieties $X$ lie in the tautological ring
by the localization formula for the 
virtual class (see [GrP1]).   We do not know any  
example defined over $\bar{\mathbb Q}$ 
of a Gromov-Witten class for which the push-forward is {\em not} 
tautological.

The moduli spaces of Hurwitz covers of $\proj^1$
also define natural classes on the moduli space of curves. 
Let $g\geq 0$.
Let $\mu^1, \ldots, \mu^m$ be $m$ partitions of equal size $d$ satisfying
$$2g-2+2d = \sum_{i=1}^m \Big( d- \len(\mu^i) \Big),$$
where $\ell(\mu^i)$ denotes the length of the partition $\mu^i$.
The moduli space of Hurwitz covers, $$H_{g}(\mu^1,\ldots,\mu^m)$$
parameterizes morphisms,
$$f: C \rarr \proj^1,$$
where $C$ is
a complete, connected, nonsingular curve 
with marked profiles $\mu^1, \ldots, \mu^m$ over $m$ ordered points of the
target (and no ramifications elsewhere).  Two Hurwitz covers,
$$C\stackrel{f}{\rarr} \proj^1, \ \ C' \stackrel{f'}{\rarr} \proj^1,$$
are 
isomorphic
if there exist isomorphisms,
$$\alpha:C \rarr C', \ \ \beta: \proj^1 \rarr \proj^1,$$
which commute with $f,f'$ and respect all the markings.
The moduli space of Hurwitz covers is a dense open set of 
the compact 
moduli space of admissible covers [HM],
$$H_{g}(\mu^1,\ldots,\mu^m) \subset  
\overline{H}_{g}(\mu^1,\ldots,\mu^m).$$
Let $\rho$ denote the map to the moduli of curves,
$$\rho: \overline{H}_{g}(\mu^1,\ldots,\mu^m) \rarr {\overline M}_{g, 
\sum_{i=1}^m \len(\mu^i)}.$$
The push-forwards of the fundamental classes,
$$\rho_*\Big( \overline{H}_{g}(\mu^1,\ldots,\mu^m)\Big) \in 
A^*(\overline{M}_{g,\sum_{i=1}^m \len(\mu^i)}),$$
define classes on the moduli of curves.

The following two questions provide classical motivation for  the results 
of our paper:
\begin{enumerate}
\item[(i)] Are the push-forwards $\rho_*\Big(  
\overline{H}_{g}(\mu^1,\ldots,\mu^m)  \Big)$ tautological?
\item[(ii)] Can the push-forwards be computed?
\end{enumerate}
We settle (i-ii) in the affirmative.

While Hurwitz covers have played a basic role in the
study of the moduli of curves,
questions (i-ii) were open even for the class of the hyperelliptic locus,
 $$\rho_*\Big( \overline{H}_{g}((2), \ldots,(2))\Big),$$
in
$\overline{M}_{g,2g+2}$.

\subsection{Stable relative maps}
\subsubsection{Overview} 
We study Hurwitz covers 
in the much richer context of stable relative maps to $\proj^1$. 
Stable relative maps combine features of stable maps and admissible 
covers.
The moduli space of stable relative maps was first
introduced by Li and Ruan in [LR]. An algebraic development
can be found in [Li].
The main results of the paper concern the relationship between
tautological classes on the moduli space of stable relative maps and
the moduli space of curves.

\subsubsection{The parameterized case}
\label{parram}
We first define stable relative maps to a parameterized $\proj^1$.
Let $g,n \geq 0$. Let $\mu^1, \ldots, \mu^m$ be $m$ partitions of
equal size $d>0$.
A stable relative map,
\begin{equation*}
  \Big[  (C, p_1,\ldots,p_n, Q_1,\ldots,Q_m) \stackrel{f}{\rarr} (T, 
q_1,...,q_m) 
\stackrel{\epsilon}{\rarr} \proj^1 \Big] \ \in  \ 
\overline{M}^{\dagger}_{g,n}(\mu^1, \ldots, \mu^m),
\end{equation*}
consists of the following data:
\begin{enumerate}
\item[(i)]  $T$ is a complete, connected genus 0 nodal curve with $m$ distinct
       nonsingular markings $q_1, \ldots, q_m$. 
\item[(ii)] The structure map
      $\epsilon: T \rarr \proj^1$ 
       restricts to 
       an isomorphism on a unique component $P\subset T$ and contracts 
all other components of $T$.
\item[(iii)]
       All
       the extreme components of  the tree $T$ not equal to $P$ carry at 
least 1 of the $m$ markings.
\item[(iv)] $C$ is a complete, connected genus $g$ nodal curve with 
$$n+ \sum_{i=1}^m \len(\mu^i)$$
distinct nonsingular markings 
$$\{p_1,\ldots,p_n\} \ \cup \ \bigcup_{i=1}^m Q_i,$$
where $|Q_i|= \len(\mu^i)$. 
\item[(v)] The morphism $f$ satisfies the following basic conditions:
\begin{enumerate}
           \item[(a)]   $f$ satisfies admissible cover conditions over the
                nodes of $T$: matching branchings with no markings or
                contracted components of $C$ lying over the nodes of $T$, 
\item[(b)]  $f$ has profile $\mu^i$ over $q_i$, 
\item[(c)] $Q_i$ is a complete marking of the fiber of $f$ over $q_i$ with  
$\len(\mu^i)$ distinct
                points.
\end{enumerate}
 \item[(vi)] The data has a {\em finite} automorphism group. 
The automorphism group is determined by
curve automorphisms,
$$ \alpha: C \rarr C,  \  \ \beta: T \rarr T,$$
which respect the markings and commute with $f,\epsilon$, and the identity 
map on $\proj^1$.
\end{enumerate}

The data
$C\stackrel{f}{\rarr} T \stackrel{\epsilon}{\rarr} \proj^1$
and  
$ 
C'\stackrel{f'}{\rarr} T' \stackrel{\epsilon'}{\rarr} \proj^1$
are isomorphic if there exist isomorphisms
$$ \alpha: C \rarr C',  \  \ \beta: T \rarr T',$$
which respect all the markings and commute with $f,f',\epsilon,\epsilon',$ 
and the identity map on $\proj^1$.
Condition (vi) may be interpreted as a stability condition.

The superscript $\dagger$ in
the notation for the moduli space of stable relative maps,
$$\overline{M}^\dagger_{g,n}(\mu^1, \ldots, \mu^m),$$ 
indicates the target $\proj^1$ is
parameterized.
The moduli space, a Deligne-Mumford stack, admits
several canonical structures:
\begin{enumerate}
\item[(i)] a virtual fundamental class,
$$[\overline{M}^\dagger_{g,n}(\mu^1, \ldots, \mu^m)]^{vir} \in 
A_{e}(\overline{M}^\dagger_{g,n}(\mu^1, \ldots, \mu^m)),$$
in the expected dimension,
$$e=2g-2+2d+n+\sum_{i=1}^m \Big( 1+\len(\mu^i) - d \Big),$$
\item[(ii)] evaluation maps,
$$ \text{ev}: \overline{M}^\dagger_{g,n}(\mu^1, \ldots, \mu^m) \rarr 
\proj^1,$$
determined by the markings $p_1, \ldots, p_n$, and $q_1, \ldots, q_m$,
\item[(iii)] cotangent line classes $\psi\in A^1( 
\overline{M}^\dagger_{g,n}(\mu^1, ..., \mu^m))$
determined by the markings $p_1, \ldots,p_n$, and $q_1,\ldots,q_m$,
\item[(iv)]  a map to the moduli of curves (via the domain),
$$\rho: \overline{M}^\dagger_{g,n}(\mu^1, \ldots, \mu^m) \rarr 
\overline{M}_{g,n+ \sum_{i=1}^{m} \len(\mu^i)   },$$
in case $2g-2+n+ \sum_{i=1}^{m} \len(\mu^i)>0$,
\item[(v)]  and a map to the Fulton-MacPherson parameter space of points 
on $\proj^1$ (via the
target),
$$\rho_0: \overline{M}^\dagger_{g,n}(\mu^1, \ldots, \mu^m) \rarr 
\proj^1[m].$$
\end{enumerate}

A {\em relative Gromov-Witten class} on the moduli space of stable relative
maps,
$$\omega =\prod_{i=1}^n \text{ev}_{p_i}^*(\gamma_i)\  \psi_{p_i}^{k_i} \ 
\cup \ 
\prod_{j=1}^m \text{ev}_{q_j}^*(\gamma'_j) \ \psi_{q_j}^{k'_j},$$
is constructed from structures (ii-iii). Here, $\gamma_i, \gamma_j' \in 
A^*(\proj^1)$.

The map to the moduli of curves (iv) may be defined to exist in
{\em all} cases relevant to our study. We will only consider moduli spaces
of stable relative maps for which $m>0$. 
Since $d$ is positive, $\len(\mu^i)>0$. Therefore, the inequality,
$$2g-2+n+ \sum_{i=1}^{m} \len(\mu^i)>0,$$
can be violated only if $g=0$. The map to the moduli of curves
is defined in the unstable genus 0 cases by viewing $\overline{M}_{0,1}$ and
$\overline{M}_{0,2}$ as points.

\subsubsection{The unparameterized case}
\label{unparram}
We will also require stable relative maps to an unparameterized $\proj^1$.
Let $g,n \geq 0$. Let $\mu^1, \ldots, \mu^m$ be $m\geq 2$ partitions of
equal size $d>0$.
A stable relative map,
\begin{equation*}
  \Big[  (C, p_1,\ldots,p_n, Q_1,\ldots,Q_m) \stackrel{f}{\rarr} (T, 
q_1,...,q_m) 
 \Big] \ \in  \ \overline{M}_{g,n}(\mu^1, \ldots, \mu^m),
\end{equation*}
consists of the following data:
\begin{enumerate}
\item[(i)]  $T$ is a complete, connected genus 0 nodal curve with $m$ distinct
       nonsingular markings $q_1, \ldots, q_m$. 
\item[(ii)]
       All
       the extreme components of  the tree $T$ carry at least 1 of the $m$ 
markings.
\item[(iii)] $C$ is a complete, connected genus $g$ nodal curve with 
$$n+ \sum_{i=1}^m \len(\mu^i)$$
distinct nonsingular markings 
$$\{p_1,\ldots,p_n\} \ \cup \ \bigcup_{i=1}^m Q_i,$$
where $|Q_i|= \len(\mu^i)$. 
\item[(iv)] The morphism $f$ satisfies (a), (b), and (c) of the 
parameterized definition.
\item[(v)] The data has a {\em finite} automorphism group. 
The automorphism group is determined by
curve automorphisms,
$$ \alpha: C \rarr C,  \  \ \beta: T \rarr T,$$
which respect the markings and commute with $f$.
\end{enumerate}

We will consider the unparameterized case 
only when $m\geq 2$.
The moduli space $\overline{M}_{g,n}(\mu^1, \ldots, \mu^m)$
admits the following structures:
\begin{enumerate}
\item[(i)] a virtual fundamental class,$$[\overline{M}_{g,n}(\mu^1, 
\ldots, \mu^m)]^{vir} \in 
A_{e}(\overline{M}_{g,n}(\mu^1,  \ldots, \mu^m)),$$
in the expected dimension,
$$e=2g-5+2d+n+\sum_{i=1}^m \Big( 1+\len(\mu^i) - d \Big),$$
\item[(ii)] cotangent line classes $\psi\in A^1( \overline{M}_{g,n}(\mu^1, 
..., \mu^m))$
determined by the markings $p_1, \ldots,p_n$, and $q_1,\ldots,q_m$,
\item[(iii)]  a map to the moduli of curves (via the domain),
$$\rho: \overline{M}_{g,n}(\mu^1, \ldots, \mu^m) \rarr 
\overline{M}_{g,n+\sum_{i=1}^{m} \len(\mu^i)},$$
in case $2g-2+n+ \sum_{i=1}^{m} \len(\mu^i)>0$,
\item[(iv)] and a map to the moduli of genus 0 pointed curves (via the 
target),
$$\rho_0: \overline{M}_{g,n}(\mu^1, \ldots, \mu^m) \rarr  
\overline{M}_{0,m},$$
in case $m\geq 3$.
\end{enumerate}

The moduli spaces of admissible covers arise as special cases of the 
moduli spaces of
stable relative maps:
$$\overline{M}_{g,0}(\mu^1, \ldots, \mu^m) =
\overline{H}_{g}(\mu^1,\ldots,\mu^m),$$
in case 
$$2g-2+2d = \sum_{i=1}^m \Big( d- \len(\mu^i) \Big).$$
Here, the virtual class of the space of stable relative maps (i) is equal to 
the fundamental
class of the space of admissible covers.

A {\em relative Gromov-Witten class} on the moduli space of stable relative
maps,
$$\omega =\prod_{i=1}^n  \psi_{p_i}^{k_i} \ \cup \ 
\prod_{j=1}^m \psi_{q_j}^{k'_j},$$
is constructed from the cotangent lines classes (ii). 
The evaluation maps are lost in the unparameterized case.
The map to the moduli of curves (iii) is defined
in all relevant cases as before.

\subsubsection{Results}
We prove results which show the compatibility of Gromov-Witten classes
on the moduli of stable relative maps and tautological classes on
the moduli of curves in both the
parameterized and unparameterized cases.

\begin{tm} 
\label{one}
Relative Gromov-Witten classes classes push-forward to tautological 
classes:
\begin{enumerate}
\item[(i)]
For every relative Gromov-Witten class $\omega$ on
 $\overline{M}^\dagger_{g,n}(\mu^1, \ldots, \mu^m)$,
$$\rho_*\Big( \omega \cap  [\overline{M}^\dagger_{g,n}(\mu^1, \ldots, 
\mu^m)]^{vir}\Big)$$
is a tautological class on the target moduli space of curves.
\item[(ii)]
For every relative Gromov-Witten class $\omega$ on
 $\overline{M}_{g,n}(\mu^1, \ldots, \mu^m)$,
$$\rho_*\Big( \omega \cap  [\overline{M}_{g,n}(\mu^1, \ldots, 
\mu^m)]^{vir}\Big)$$
is a tautological class on the target moduli space of curves.
\end{enumerate}
\end{tm}

Theorem \ref{one} is proven by studying relations obtained by virtual 
localization.
The proof for both parts  is {\em constructive}:  the
push-forwards of relative Gromov-Witten classes are recursively
calculated in the tautological ring in terms of the
standard $\psi$, $\kappa$, and $\lambda$ classes.
The proof is given in Sections \ref{argg}-\ref{matt} of the paper.

%Stable relative maps to $\proj^1$ may be defined in the
%degree 0 case with partitions $\mu^i= \emptyset$. Theorem \ref{one}
%in the degree 0 case follows from a direct examination of
%the virtual classes of these simple spaces. 

\subsection{Consequences}
Theorem \ref{one}  has several consequences for the
geometry of the moduli space of curves. 

\begin{pr}
\label{j}
The moduli of Hurwitz covers yields tautological classes,
$$\rho_*\Big( {\overline{H}}_g(\mu^1, \ldots, \mu^m)\Big) 
\in R^*(\overline{M}_{g, \sum_{i=1}^m 
\len(\mu^i)}).$$
\end{pr}

Proposition \ref{j} follows from part (ii) of Theorem \ref{one} in case
$\omega$ is the identity class and the space of stable relative
maps is specialized to the space of admissible covers. Moreover, the 
push-forwards
are effectively determined. 

An additive set of generators for $R^*(\overline{M}_{g,n})$ is indexed
by strata dual graphs decorated with $\psi$ and $\kappa$ classes on the
nodes, see [GrP2].
We review the basic push-pull method 
for constructing relations in $R^*(\overline{M}_{g,n})$ introduced in [P1], 
[BP].
Consider the two maps to the moduli of curves
defined for the moduli space of admissible covers:
$$\overline{M}_{0,m} \stackrel{\rho_0}
{\longleftarrow}
\overline{H}_g(\mu^1,\ldots,\mu^m) \stackrel{\rho}{\longrightarrow} 
\overline{M}_{g,\sum_{i=1}^n \len(\mu^i)}.$$
Let $r$ denote a relation among boundary strata in $\overline{M}_{0,m}$. 
Then,
\begin{equation}
\label{rell}
\rho_* \rho_0^{*}(r),
\end{equation}
defines a relation in $\overline{M}_{g, \sum_{i=1}^m \len(\mu^i)}$.

Relation (\ref{rell}), however, consists of push-forwards of 
fundamental classes of auxiliary moduli spaces of admissible covers.
The push-pull method, along with a case-by-case 
analysis of the push-forwards,
yields the basic relations
in genus 1 and 2 among descendent stratum classes, see [P1], [BP], [Ge1], 
[Ge2].
Proposition \ref{j} now guarantees all push-forwards arising
in (\ref{rell}) lie in $R^*(\overline{M}_{g,n})$.
The push-pull method together with Proposition \ref{j} provides
a rich source of tautological relations in the moduli space of curves.

\medskip
\noindent {\bf Speculation 1.}  {\em All relations in the tautological 
ring are obtained
via the push-pull method 
and Proposition \ref{j}.}
\medskip

A further study of the push-pull relations was undertaken in [I].  
In the absence of Proposition \ref{j}, only the principal terms  (the {\em 
symbols}) of the
relations could be studied. The main result of [I] is
of interest here:
\begin{enumerate}
\item[(i)] Ionel's vanishing,
$$\prod_{i=1}^n \psi_i^{e_i} \prod_{j\geq 0} \kappa_j^{f_j} = 0 
\in R^*(M_{g,n}), 
\ \ \text{for} \ 
\sum_{i=1}^n e_i + \sum_{j\geq 0} j f_j \geq  g + \delta_{0g} - \delta_{0n},$$
 can be obtained from the symbols of push-pull relations.
\end{enumerate}
The above vanishing generalizes two well-known vanishings:
\begin{enumerate}
\item[(ii)] 
Getzler's (conjectured) vanishing,
$$\prod_{i=1}^n \psi_i^{e_i} = 0 \in R^*(M_{g,n}), \ \ \text{for} \ 
\sum_{i=1}^n e_i \geq g + \delta_{0g},$$
\item[(iii)]
 Looijenga's vanishing,
$$\prod_{j\geq 0} \kappa_j^{f_j} = 0 \in R^*(M_g), \ \ \text{for} \ 
\sum_{j\geq 0} j f_j \geq g-1.$$
\end{enumerate}

The tautological ring, $R^*(M_g) \subset R^*(\overline{M}_g),$
for the open moduli space $M_g$, is
defined to be the image of $R^*(\overline{M}_g)$ via the natural map,
$$R^*(\overline{M}_g) \subset A^*(\overline{M}_{g}) \rarr A^*(M_g) .$$
Let $\partial \overline{M}_g= \overline{M}_g \setminus M_g$ denote
the boundary of the moduli space of curves.
The tautological classes of the boundary,
 $R^*(\partial\overline{M}_g) \subset A^*(\partial \overline{M}_g)$,
are defined by the push-forwards of all tautological classes on the
normalized boundary divisors via the gluing morphisms.
There is a restriction sequence:
$$R^*(\partial \overline{M}_g) \rarr R^*(\overline{M}_g) \rarr R^*(M_g) 
\rarr 0,$$
for which the exactness in the middle is {\em unknown}.
An identical discussion holds for $\overline{M}_{g,n}$.

Proposition \ref{j} combined with the symbol analysis of [I]
directly yields a much stronger version of Ionel's
vanishings.

\begin{pr} \label{jj}
The class $  \prod_{i=1}^n \psi_i^{e_i} \prod_{j\geq 0} \kappa_j^{f_j}  $ 
lies in the image  
$$Im \ R^*(\partial \overline{M}_{g,n}) \subset R^*(\overline{M}_{g,n}),$$
for $ \sum_{i=1}^n e_i + \sum_{j\geq 0} j f_j \geq  g
+ \delta_{0g} - \delta_{0n}. $
\end{pr}

Proposition \ref{jj} has consequences beyond Ionel's original 
vanishing.
A uniform Gorenstein conjecture was advanced in [FP], [P2] for the 
tautological
rings of $\overline{M}_{g,n}$ and the quotients corresponding to the 
moduli
spaces of compact type curves $M_{g,n}^c$
and curves with rational tails $M_{g,n}^{rt}$.
An inductive argument using Proposition \ref{jj} yields the 
following
result.

\begin{pr} \label{jjjj}
The socle and vanishing claims of the Gorenstein conjectures hold for
$\overline{M}_{g,n}$, $M_{g,n}^c$, and $M_{g,n}^{rt}$.
\end{pr}

Proposition \ref{jjjj} is proven in Section \ref{vs}. 
A different approach to Proposition
\ref{jjjj} has been pursued recently by Graber and Vakil [GrV1], [GrV2].

Another consequence of Proposition \ref{jj} is the exactness of the
restriction sequence,
\begin{equation}
\label{rigg}
R^*(\partial \overline{M}_g) \rarr R^*(\overline{M}_g) \rarr R^*(M_g) 
\rarr 0,
\end{equation}
in degrees greater than $g-2$. The result motivates the following conjecture.

\medskip
\noindent {\bf Conjecture 2.}  {\em The restriction sequence (\ref{rigg}) is 
exact in all degrees.}

\medskip\noindent
More generally, we conjecture the exactness of the tautological sequences
associated to the compact type and rational tail moduli spaces, 
$$ R^*(\overline{M}_{g,n} \setminus M^c_{g,n}) 
\rarr R^*(\overline{M}_{g,n}) \rarr R^*(M^c_{g,n}) 
\rarr 0,$$
$$ R^*(\overline{M}_{g,n} \setminus M^{rt}_{g,n}) 
\rarr R^*(\overline{M}_{g,n}) \rarr R^*(M^{rt}_{g,n}) 
\rarr 0,$$
see [FP], 
[P2].

%If the push-pull method suffices to generate all relations among the 
%$\kappa$
%classes in $R^*(M_g)$, then Conjecture 2 follows from Proposition \ref{j}.

Finally, Proposition \ref{jj} 
is the {\em required} form of Getzler's 
vanishing
for applications to Gromov-Witten theory. The main consequence is the
following reconstruction result proved in Section \ref{gw}.

\begin{pr}\label{jjjjj}
Let $X$ be a nonsingular projective variety. All descendent Gromov-Witten 
invariants
of $X$ can be reconstructed from the restricted invariants
$$\langle \ \tau_{e_1}(\gamma_1)\cdots \tau_{e_n}(\gamma_n) \prod_{j\geq 0} 
\kappa_j^{f_j}
\ \rangle_{g,n,\beta}^X,$$
where $\sum_{i=1}^n e_i + \sum_{j\geq 0} j f_j < g +\delta_{0g}$.
\end{pr}

We conjecture a stronger result based on a new conjectured
generation statement for $R^*(M_{g,n})$ and 
the conjectured exactness of tautological restriction sequences.
The conjectural framework is discussed in Section \ref{lll}.

\medskip
\noindent {\bf Conjecture 3.} {\em All descendent Gromov-Witten 
invariants
of $X$ can be reconstructed from the restricted invariants
$$\langle \ \tau_{e_1}(\gamma_1)\cdots \tau_{e_n}(\gamma_n) 
\ \rangle_{g,n,\beta}^X,$$
where $\sum_{i=1}^n e_i  < g + \delta_{0g} $.}

\subsection{Acknowledgements}

We thank T. Graber, A. Okounkov, Y. Ruan, and R. Vakil for discussions of the 
relative Gromov-Witten theory of $\proj^1$. 
The research reported here was pursued during a visit by C.F. to
Princeton University in the spring of 2003.
R.P. was partially supported by DMS-0236984
and fellowships from the Sloan and Packard foundations.

\section{Localization}
\subsection{Overview}
The localization formula for the virtual class of the moduli space
$\overline{M}^\dagger_{g,n}(\mu^1, \ldots, \mu^m)$ is
required for our proof of Theorem \ref{one}. The
localization formula can be obtained from the algebraic construction of
the virtual class [Li] together with the virtual localization formula
[GrP1]. A derivation can be found in 
[GrV2].

\subsection{Disconnected domains}

\subsubsection{Notation} We will require moduli spaces of stable relative maps with
{\em disconnected} domains in both the parameterized and
unparameterized cases. The definitions follow the connected case with
minor variations.

We first introduce notation for the disconnected case.
Let $\GG$ be an ordered set of genera,
 $$\GG=(g_1, \ldots, g_c).$$
Let $\NN$ be an ordered set partition,
$$\NN=(n_1,\ldots,n_c), \ \ \bigcup_{i=1}^c n_i =  \underline{n},$$
where $\underline{n}= \{1, \ldots,n\}.$
The data $\GG$ and $\NN$ describe the  genera and marking
distributions of a disconnected domain with $c$ ordered components.

The degree distribution on a disconnected domain is described by
an ordered partition $\DD$,
$$\DD= (d_1,\ldots, d_c), \ \ d_i >0, \ \ \sum_i d_i = d.$$ 
A partition $\mumu$ of $d$ of {\em type} $\DD$ is an ordered
set of partitions,
$$\mumu= (\mu[1], \ldots, \mu[c]),$$
where $\mu[i]$ is a partition of $d_i$.
Ordered partitions of type $\DD$ describe ramification conditions
on the disconnected domain.

\subsubsection{Moduli spaces}
We first define stable relative maps to a parameterized $\proj^1$
with disconnected domains.
Let $\GG$, $\NN$, and $\DD$ be as defined above. 
Let $\mumu^1, \ldots, \mumu^m$ be $m$ partitions of
$d$ of type $\DD$.
A stable relative map,
\begin{equation*}
  \Big[  (C, p_1,\ldots,p_n, Q_1,\ldots,Q_m) \stackrel{f}{\rarr} (T, 
q_1,...,q_m) 
\stackrel{\epsilon}{\rarr} \proj^1 \Big] \ \in  \ 
\overline{M}^{\dagger}_{\GG,\NN}(\mumu^1, \ldots, \mumu^m),
\end{equation*}
consists of the data (i-vi) of Section \ref{parram} with 
(iv) replaced by:
\begin{enumerate}
\item[(iv)] $C$ is a complete, disconnected
nodal curve with $c$ ordered components carrying
$$n+ \sum_{i=1}^m \len(\mumu^i)$$
distinct nonsingular markings 
$$\{p_1,\ldots,p_n\} \ \cup \ \bigcup_{i=1}^m Q_i,$$
where $|Q_i|= \len(\mumu^i)$.
The genera of the components are determined by $\GG$. The
markings are distributed by the data of $\NN$ and $\mumu^i$.
\end{enumerate}

The data
$C\stackrel{f}{\rarr} T \stackrel{\epsilon}{\rarr} \proj^1$
and  
$ 
C'\stackrel{f'}{\rarr} T' \stackrel{\epsilon'}{\rarr} \proj^1$
are isomorphic if there exist isomorphisms
$$ \alpha: C \rarr C',  \  \ \beta: T \rarr T',$$
which respect all the structures and commute with $f,f',\epsilon,\epsilon',$ 
and the identity map on $\proj^1$.

The moduli space, a Deligne-Mumford stack, admits the 
canonical structures  (i-iv) of the connected case. 
Relative Gromov-Witten classes are defined as before using structures (ii-iii).
The map to the moduli
of curves (iv) via the domain  is slightly altered:
$$\rho: \overline{M}^\dagger_{\GG,\NN}(\mumu_1, \ldots, \mumu_m) \rarr 
\prod_{i=1}^c \overline{M}_{g_i,|n_i|+ \sum_{j=1}^{m} \len(\mu^j[i])   },$$
in case 
$$2g_i-2+|n_i|+ \sum_{j=1}^{m} \len(\mu^j[i])>0,$$
for all $i$.
Since $m>0$ and $d_i>0$, the map to the moduli of curves
is defined in all relevant cases by viewing $\overline{M}_{0,1}$ and
$\overline{M}_{0,2}$ as points.

The definition of the unparameterized
moduli space ${\overline{M}}_{\GG,\NN}(\mumu^1,
\ldots,\mumu^m)$
is obtained similarly. Part (iii) of the definition of 
Section \ref{unparram} is modified to allow
disconnected domains. Also,
the canonical map to the moduli of curves (iii) is replaced by a map to
a product of moduli spaces.

Of course, when $c=1$, the moduli spaces of stable relative maps
with connected domains are recovered.

\subsubsection{Theorem \ref{one} revisited}
We will prove a stronger form of Theorem 1 including all the
moduli spaces of stable relative maps with disconnected domains.

\begin{tm} 
\label{two}
Relative Gromov-Witten classes push-forward to tautological 
classes:
\begin{enumerate}
\item[(i)]
For every relative Gromov-Witten class $\omega$ on
 $\overline{M}^\dagger_{\GG,\NN}(\mumu^1, \ldots, \mumu^m)$,
$$\rho_*\Big( \omega \cap  [\overline{M}^\dagger_{\GG,\NN}(\mumu^1, \ldots, 
\mumu^m)]^{vir}\Big)$$
is a tautological class on the target.
\item[(ii)]
For every relative Gromov-Witten class $\omega$ on
 $\overline{M}_{\GG,\NN}(\mumu^1, \ldots, \mumu^m)$,
$$\rho_*\Big( \omega \cap  [\overline{M}_{\GG,\NN}(\mumu^1, \ldots, 
\mumu^m)]^{vir}\Big)$$
is a tautological class on the target.
\end{enumerate}
\end{tm}

\subsection{Torus actions, fixed points, and the localization formula}

\subsubsection{The torus action}
The equivariant Chow ring of $\com^*$ is freely generated by
the first Chern class $t$ of the standard representation,
$$A^*_{\com^*} ([point]) = {\Q}[t].$$

Let $\com^*$ act diagonally on a two-dimensional vector
space $V$ via the trivial and
standard representations,
\begin{equation}
\label{repp}
\xi\cdot (v_1,v_2) = ( v_1, 
\xi \cdot v_2).
\end{equation}
Let $\proj^1=\proj(V)$.
Let $0,\infty $ be the fixed points $[1,0], [0,1]$ of the corresponding
$\com^*$-action on $\proj(V)$.

An equivariant lifting  of $\com^*$ to a line bundle $L$
over 
$\proj(V)$ is uniquely determined by
the fiber
representations at the fixed points 
$$L_{0}, \ \  L_{\infty}.$$
The canonical lifting of $\com^*$ to the
tangent bundle $T_\proj$ has weights $[t,-t]$.

The representation (\ref{repp}) canonically
induces a $\com^*$-action on the moduli space 
$\overline{M}_{\GG,\NN}(\mumu^1, \ldots, \mumu^m)$  by translation of
the map:
$$\xi \cdot[C\stackrel{f}{\rarr} T \stackrel{\epsilon}{\rarr} \proj^1]
= [C\stackrel{f}{\rarr} T \stackrel{\xi \cdot \epsilon}{\rarr} \proj^1].$$

The canonical structures of the moduli space of stable relative maps
are
compatible with the induced $\com^*$-action.
The virtual fundamental class canonically lifts to
equivariant Chow theory,
$$[\overline{M}^\dagger_{\GG,\NN}(\mumu^1, \ldots, \mumu^m)]^{vir} \in 
A^{\com^*}_{e}(\overline{M}^\dagger_{g,n}(\mu^1, \ldots, \mu^m)).$$
The $\com^*$-action on the moduli space canonically lifts to the cotangent
line bundles and is equivariant with respect to the evaluation maps.
Therefore, equivariant relative Gromov-Witten classes,
$$\omega =\prod_{i=1}^n \text{ev}_{p_i}^*(\gamma_i)\  \psi_{p_i}^{k_i} \ 
\cup \ 
\prod_{j=1}^m \text{ev}_{q_j}^*(\gamma'_j) \ \psi_{q_j}^{k'_j},$$
are well-defined for $\gamma_i, \gamma_j' \in 
A_{\com^*}^*(\proj^1)$. 
The $\com^*$-action on the moduli space
is equivariant (via $\rho$) with respect to the trivial action
on the moduli space of curves.

We will require a localization formula for the equivariant push-forward:
$$\rho_*\Big(\omega \cap [\overline{M}^\dagger_{\GG,\NN}(\mumu^1, \ldots, 
\mumu^m)]^{vir}\Big).$$

\subsubsection{Graph notation}
Let the data $$\GG=(g_1, \ldots,g_c), \ \ \NN=(n_1, \ldots,n_c),
\ \ \DD=(d_1, \ldots,d_c)$$ describe
(possibly) disconnected domains with $c$ components.
The torus fixed loci of $\overline{M}^\dagger_{\GG,\NN}(\mumu^1, \ldots, 
\mumu^m)$ are indexed by {\em localization graphs},
$$\Gamma =(V,E,N, \gamma,\pi,\delta,(R^1,\ldots, R^m)),$$
defined by the following conditions:
\begin{enumerate}
\item[(i)] $V$ is the vertex set,
\item[(ii)] $\gamma: V \rarr {\mathbb{Z}}_{\geq 0}$ is a genus assignment,
\item[(iii)] $\pi: V \rarr \{0,\infty\}$ is a function,
\item[(iv)] $E$ is the edge set, 
\begin{enumerate}
\item[(a)]
If an edge $e$ connects $v,v'\in V$, then
$\pi(v)\neq \pi(v')$, \\
in particular, there are no self edges, 
\item[(b)]
$\Gamma$ has $c$ {\em ordered} connected components 
$$\Gamma_1, \ldots, \Gamma_c,$$
with vertex and edge sets $(V_i, E_i)$ respectively,
\end{enumerate}
\item[(v)] $\delta: E \rarr {\mathbb{Z}}_{>0}$ is a degree
assignment,
\item[(vi)] $N=\{1, \ldots,n\}$ is a set of vertex markings for which
each subset $n_i$ is incident to $V_i$,
\item[(vii)] $g_i= \sum_{v\in V_i} \gamma(v) + h^1(\Gamma_i)$,
\item[(viii)] $d_i= \sum_{e\in E_i} \delta(e).$
\item[(ix)] $R^j$ is a {\em refinement} of $\mumu^j$ consisting of:
\begin{enumerate} 
\item[(a)]   a choice of {\em side} $s^j \in \{0, \infty\}$,
\item[(b)]  a distribution of the parts of $\mumu^j$ to the
vertices of $\pi^{-1}(s^j)$ satisfying two properties:
\begin{enumerate}
\item[$\bullet$] the parts of $\mumu^j[i]$ are
distributed to the vertices of $\Gamma_i$,
\item[$\bullet$] the sum of the parts distributed to $v$ equals the sum
                  of the degrees of edges incident to $v$.
\end{enumerate}
\end{enumerate}
\end{enumerate}

\subsubsection{Torus fixed points}
Let the moduli point,
\begin{equation} 
\label{fxxx}
[C\stackrel{f}{\rarr} T 
\stackrel{\epsilon}{\rarr} \proj^1] 
\in \overline{M}^\dagger_{\GG,\NN}(\mumu^1, \ldots, \mumu^m),
\end{equation}
be fixed by the torus action. 
All marked points, nodes, contracted components, and ramification
points of $C$ must lie over the torus fixed set $\{ 0,\infty\}$ of $\proj^1$.   
Each irreducible component $D\subset C$ dominant onto $\proj^1$  
must be a Galois cover with full ramification 
over the two fixed points $\{0,\infty\}$.

We associate a localization
graph, 
$$\Gamma=(V,E,N, \gamma,\pi,\delta,(R^1,\ldots, R^m)),$$ to the torus fixed point (\ref{fxxx})
by the  following construction:
\begin{enumerate}
\item[(i)] Let 
$V$ be the set of connected components of 
$(\epsilon \circ f)^{-1}(\{ 0,\infty\})$, 
\item[(ii)] Let 
$\gamma(v)$ be the arithmetic genus of the corresponding component 
(taken to be 0 if the component is an isolated point),
\item[(iii)] Let $\pi(v)$ be 
the fixed point in $\proj^1$ over which the corresponding component lies,
\item[(iv)] Let $E$ be 
the set of non-contracted irreducible components $D\subset C$,
\item[(v)] Let $\delta(D)$ be the degree of the Galois cover $\epsilon \circ f|_D$,
\item[(vi)] Let $N$ be the marking set.
\item[(ix)] Let $R^j$ be the refinement of $\mumu^j$ obtained
by the ramification
            conditions.
\end{enumerate}
All the conditions of a localization graph, including (vii-viii), hold by 
the definition of a stable relative map.

A stack $\overline{M}_\Gamma$ together with an action of
a finite group ${\mathbf A}_\Gamma$ is canonically constructed
from $\Gamma$ in Section \ref{nxx} below. A canonical inclusion will be defined,
$$\tau_\Gamma/{\mathbf A}_\Gamma: \overline{M}_\Gamma/ {\mathbf A}_\Gamma \rarr
\overline{M}^\dagger_{\GG,\NN}(\mumu^1, \ldots, \mumu^m).$$ 
The disjoint union,
$$\bigcup_{\Gamma}\ \overline{M}_\Gamma/{\mathbf A}_\Gamma,$$
will equal the total torus fixed set.

\subsubsection{The torus fixed locus $\overline{M}_\Gamma/{\mathbf A}_\Gamma$}
\label{nxx}
Let $\Gamma$ be a bipartite graph.
The stack $\overline{M}_\Gamma$ is defined as a product of auxiliary
moduli spaces of curves and maps. The  ${\mathbf A}_\Gamma$-action
is obtained from the automorphisms of $\Gamma$. 

\vspace{+10pt}
\noindent {\bf Case I.} The refinements $(R^1, \ldots, R^m)$
lie on {\em both} sides $0$ and $\infty$.
\vspace{+10pt}

The data of $\Gamma$ over $0$ uniquely defines a
moduli space of unparameterized stable relative maps,
$$\overline{M}_0= \overline{M}_{\GG_0,\NN_0}(\mathbf{R}_0, R_\delta),$$
where:
\begin{enumerate}
\item[(i)] $\GG_0$ is determined by the genera of the vertices $V_0= f^{-1}\{0\}$,
\item[(ii)] $\NN_0$ is determined by the markings of $V_0$,
\item[(iii)] ${\mathbf R}_0$ is the set of refinements on side $0$,
\item[(iv)] $R_\delta$ is the ramification condition determined by $\delta$.
\end{enumerate}
Let $q_0$ denote the new marking associated to $R_\delta$.

Similarly, the unparameterized moduli space,
    $$\overline{M}_\infty= \overline{M}_{\GG_\infty,\NN_\infty}(R_\delta,
\mathbf{R}_\infty),$$
is defined by the data of $\Gamma$ over $\infty$.
Let $q_\infty$ denote the new marking associated to $R_\delta$.

Let $\overline{M}_\Gamma$ be the product,
$\overline{M}_\Gamma= \overline{M}_0 \times \overline{M}_\infty.$
The moduli space $\M_\Gamma$ has a virtual class,
$$[\M_\Gamma]^{vir} = [\overline{M}_0]^{vir} \times [\overline{M}_\infty]^{vir},$$
 determined
by the product of the virtual classes of the two factors.

Over $\M_\Gamma$, there is
a canonical family of
$\com^*$-fixed stable relative maps, $$\pi_C: {\mathcal C}\rarr \M_\Gamma,$$
$$\pi_T: {\mathcal T} \rarr \M_\Gamma,$$
$${\mathcal C} 
\stackrel{f}{\rarr} {\mathcal T} \stackrel{\epsilon}{\rarr} \proj^1.$$
The canonical family is
constructed by attaching the universal families over $\M_0$ and $\M_\infty$
according to $\Gamma$. The canonical family
yields a canonical morphism
of stacks,
 $\tau_\Gamma: \M_\Gamma 
\rarr \overline{M}^\dagger_{\GG,\NN}(\mumu^1, \ldots, \mumu^m)
.$

There is a natural automorphism group ${\mathbf A}_\Gamma$ acting 
equivariantly on ${\mathcal{C}}$ and 
$\M_\Gamma$ with respect to the morphisms $\rho$ and
$\pi_C$. The group
${\mathbf A}_\Gamma$ acts via
automorphisms of the Galois covers (corresponding to the
edges) and the symmetries of the graph $\Gamma$.
The group 
${\mathbf A}_\Gamma$ is filtered by an exact sequence of groups,
$$ 1 \rarr \prod_{e\in E} {\Z}/{\delta(e)} \rarr
{\mathbf A}_\Gamma \rarr \text{Aut}(\Gamma) \rarr 1,$$
where $\text{Aut} (\Gamma)$ is the 
automorphism group of $\Gamma$: $\text{Aut}(\Gamma)$ is the
subgroup of the permutation group of the vertices and edges which
respects all the structures of $\Gamma$.
$\text{Aut}(\Gamma)$ naturally acts on $\prod_{e\in E} \Z/ \delta(e)$
and ${\mathbf A}_\Gamma$ is the semidirect product. 

The quotient stack $\M_\Gamma/{\mathbf A}_\Gamma$
is a nonsingular Deligne-Mumford stack.
The induced 
map,
$$\tau_\Gamma/ {\mathbf A}_\Gamma
 : \M_\Gamma/{\mathbf A}_\Gamma
 \rarr \overline{M}^\dagger_{\GG,\NN}(\mumu^1, \ldots, \mumu^m),
     $$
is a closed immersion of Deligne-Mumford stacks.

The multiplicity $m(\Gamma)$ and
the Euler class of the virtual normal bundle $e(N^{vir}_\Gamma)$ will
be required for the localization formula:
$$m(\Gamma) = \prod_{e\in E} \delta(e)^2,$$
$$
\frac{1}{e(N_\Gamma^{vir})} = \frac{1}{t (t-\psi_{q_{0}})}\ \frac{1}{-t(-t-\psi_{q_{\infty}})}.$$

A degenerate configuration occurs over $0$ if the following special conditions hold for
$\Gamma$:
\begin{enumerate}
\item[(i)] all vertices of $V_0$ have genus 0,
\item[(ii)] all vertices of $V_0$ have valence 1,
\item[(iii)] a single $R^j$ lies on side 0 {\em and} $R^j=R_\delta$.
\end{enumerate}
Condition (ii) implies that no markings are incident to $V_0$.

If $\Gamma$ satisfies (i-iii), then
$\overline{M}_0$ is defined to be a point.
The multiplicity formula and Euler class formulas for the
degenerate configuration are:
$$m(\Gamma) = \prod_{e\in E} \delta(e),$$ 
$$
\frac{1}{e(N_\Gamma^{vir})} = \frac{1}{t}\ \frac{1}{-t(-t-\psi_{q_{\infty}})}.$$

Similarly, a degenerate configuration may occur over $\infty$.
The treatment is identical with the roles of 0 and $\infty$ interchanged.

In fact, degenerate configurations
may occur simultaneously at 0 and $\infty$. Then, both $\overline{M}_0$
and $\overline{M}_\infty$ are defined to be points and
$$m(\Gamma)=1,$$
$$
\frac{1}{e(N_\Gamma^{vir})} = \frac{1}{t} \ \frac{1}{-t}.$$

\vspace{+10pt}
\noindent {\bf Case II.} The refinements lie only on side $0$.
\vspace{+10pt} 

The data of $\Gamma$ over $0$ defines a
moduli space of unparameterized stable relative maps,
$$\overline{M}_0= \overline{M}_{\GG_0,\NN_0}
(\mathbf{R}_0, R_\delta),$$
as in Case I.

The data of $\Gamma$ over $\infty$ determines a product of moduli spaces,
\begin{equation}
\label{ppww}
\overline{M}_{\infty} = \prod_{v\in V_\infty} \overline{M}_{\gamma(v),val(v)}.
\end{equation}
Here, the valence $val(v)$ of a vertex counts both the incident
edges and incident markings.
The unstable moduli spaces $\overline{M}_{0,1}$ and $\overline{M}_{0,2}$
arising in the product (\ref{ppww}) are viewed as points.

Let $\overline{M}_\Gamma$ be the product, $\overline{M}_0 \times \overline{M}_\infty$.
The moduli space $\M_\Gamma$ has a virtual class,
$$[\M_\Gamma]^{vir} = [\overline{M}_0]^{vir} \times [\overline{M}_\infty],$$
 determined
by the product of the virtual class of the first factor and the
ordinary fundamental class of the second factor.

There exists a canonical family over $\M_\Gamma$ and a canonical map,
$$\tau_\Gamma:
 \M_\Gamma \rarr \overline{M}^\dagger_{\GG,\NN}(\mumu^1, \ldots, \mumu^m),
     $$
equivariant with respect to an
${\mathbf A}_\Gamma$-action exactly as in Case I.
A closed immersion,
$$\tau_\Gamma/ {\mathbf A}_\Gamma
 : \M_\Gamma/{\mathbf A}_\Gamma
 \rarr \overline{M}^\dagger_{\GG,\NN}(\mumu^1, \ldots, \mumu^m),
     $$
is obtained.

We define the multiplicity and the Euler class of
the virtual normal bundle of $\M_\Gamma$ in Case II by the formula:
$$m(\Gamma) = \prod_{e\in E} \delta(e),$$
$$\frac{1}{e(N_\Gamma^{vir})} = \frac{1}{t (t-\psi_{q_{0}})} \prod_{e\in E} 
\frac{-t}{(-t)^{\delta(e)}\frac{\delta(e)!}
{\delta(e)^{\delta(e)}}} 
\prod_{v\in V_\infty}
\frac{1}{N(v)}.$$
The vertex terms, $N(v)$, are discussed below.

 A vertex $v\in V_\infty$ is {\em stable} if
$2\gamma(v)-2+ val(v) >0.$
If $v$ is stable, the moduli space $\overline{M}_{\gamma(v), val(v)}$
is a factor of $\overline{M}_\infty$. The vertex term
$N(v)$ is an  
equivariant cohomology class on the factor $\overline{M}_{\gamma(v),val(v)}$ in
the stable case.

\vspace{10pt}
$\bullet$ Let $v\in V_\infty$ be a stable vertex.
Let $e_1, \ldots, e_l$ denote the distinct edges incident to $v$ in bijective
correspondence to a subset of the local markings of the moduli
space $\overline{M}_{\gamma(v), val(v)}$.
Let $\psi_i$ denote the cotangent line of the marking at $v$ corresponding to $e_i$,
and let $\lambda_j$ denote the Chern classes of the Hodge bundle. Then, 
$$\frac{1}{{N}(v)} = 
\frac{1}{ -t} \cdot
\prod_{i=1}^l \frac{1}
{-\frac{t}{\delta(e_i)} - \psi_i} 
\cdot 
 \sum_{j=0}^{\gamma(v)} (-1)^j \lambda_j (-t)^{\gamma(v)-j}.$$

If $v\in V_\infty$ is an unstable vertex, then $\gamma(v)=0$ and $val(v) \leq 2$.
There are three unstable cases: 
two with valence 2 and one with valence 1.

\vspace{10pt}
\noindent $\bullet$ Let $v\in V_\infty$ be an {\em unmarked} vertex with
$\gamma(v)=0$ and $val(v)=2$. Let $e_1$ and $e_2$ be the two
incident edges. Then,
$$\frac{1}{{N}(v)} 
= \frac{1}{-t} \cdot
\frac{1}{-\frac{t}{\delta(e_1)} 
- \frac{t}{\delta(e_2)}}.$$

\noindent $\bullet$ Let $v\in V_\infty$ be a $1$-marked vertex with
$\gamma(v)=0$ and $val(v)=2$. Then,
$$\frac{1}{{N}(v)} = \frac{1} {-t},$$
there are no contributing factors.

\noindent $\bullet$ Let $v\in V_\infty$ be an {\em unmarked} vertex with
$\gamma(v)=0$ and $val(v)=1$. Let $e$ be the unique
incident edge. Then,
$$\frac{1}{{N}(v)} = 
\frac{1} {-t} \cdot
\frac{-t}{\delta(e)}.$$

A degenerate configuration may arise over $0$ in Case II.
The treatment exactly follows the discussion in Case I.

\vspace{+10pt}
\noindent {\bf Case III.} The refinements lie only on side $\infty$
\vspace{+10pt}

The treatment of Case III  exactly follows the discussion of Case II with the
roles of $0$ and $\infty$ interchanged (and $t$ replacing $-t$ in all the
formulas).

\subsubsection{The localization formula}

The localization formula for the virtual class of
the moduli space $\overline{M}^\dagger_{\GG,\NN}(\mumu^1, \ldots, \mumu^m)$ is:
\begin{equation}
\label{llwww}
[\overline{M}^\dagger_{\GG,\NN}(\mumu^1, \ldots, \mumu^m)]^{vir} =
\sum_\Gamma \frac{m(\Gamma)}{|{\mathbf A}_\Gamma|} \ \tau_{\Gamma*}
\Big( \frac{
[\M_\Gamma]^{vir}}
{e(N_\Gamma^{vir})}\Big),
\end{equation}
in localized equivariant Chow theory.

\subsubsection{First application}
The localization formula  
immediately yields the 
following implication.

\begin{lm}
\label{imppp}
Theorem 2 part (i) is a consequence of Theorem 2 part (ii).
\end{lm}

\bpf
Let $\omega$ be the canonical equivariant lift of a
relative Gromov-Witten class on the
parameterized moduli space $\overline{M}^\dagger_{\GG,\NN}(\mumu^1, \ldots, 
\mumu^m)$. By the localization formula (\ref{llwww}),
\begin{equation}
\label{ffkk}
\rho_*\Big(\omega \cap [\overline{M}^\dagger_{\GG,\NN}(\mumu^1, \ldots, 
\mumu^m)]^{vir}\Big) =
\sum_\Gamma \frac{m(\Gamma)}{|{\mathbf A}_\Gamma|} \ \rho_{\Gamma*}\Big( \frac{\tau_\Gamma^*(\omega)} {e(N_\Gamma^{vir})}  \cap 
[\M_\Gamma]^{vir} \Big),
\end{equation}
where $\rho_\Gamma= \rho \circ \tau_\Gamma.$

We now analyze the $\Gamma$ term on the right side  of (\ref{ffkk}). 
The space $\overline{M}_\Gamma$ is a product of unparameterized
moduli spaces of stable relative maps and moduli spaces of stable curves.
The class $\tau_\Gamma^*(\omega)$ is composed of Gromov-Witten and
tautological classes on these factor spaces (together with powers of $t$).
Similarly, the expansion of $\frac{1}{e(N^{vir}_\Gamma)}$ is composed of
Gromov-Witten and tautological classes (together with powers of $t$).
Hence, by Theorem 2 part (ii), the equivariant term, 
$$\rho_{\Gamma*}\Big( \frac{\tau_\Gamma^*(\omega)} {e(N_\Gamma^{vir})}  \cap 
[\M_\Gamma]^{vir} \Big),$$
is a series in $t$ with coefficients in the tautological ring of the
target of $\rho$.

The non-equivariant limit of the $\rho$ push-forward,
$$\rho_*\Big(\omega \cap [\overline{M}^\dagger_{\GG,\NN}(\mumu^1, \ldots, 
\mumu^m)]^{vir}\Big),$$
is obtained from the $t^0$ coefficient of the right side of (\ref{ffkk}).
\epf

\section{Theorem 2 part ($ii$)}
\label{argg}

\subsection{Overview}
We obtain basic relations constraining the $\rho$ push-forward of relative 
Gromov-Witten
classes on the unparameterized spaces 
from the localization formula on parameterized spaces. The relations
are proven to recursively determine all the $\rho$ push-forwards in terms
of tautological classes on the target of $\rho$.

For the proof of Theorem 2 part (ii), we will require relations for 
the disconnected case. However, for ease of presentation, we first
discuss the connected case. The disconnected case follows with minor
modifications.

\subsection{Basic relations I: The connected case}

\subsubsection{The set $\Pi(d,n,k)$}
\label{pidnk0}
Let $d$ and $n$ be integers satisfying $d\ge n>0$. Let $k\ge0$ be an integer.
A partially ordered partition,
 $$\bar{\alpha}=(\alpha, \alpha'),$$ of degree $d$ and order $n$
consists of the following data:
\begin{enumerate}
\item[(i)] an ordered
partition with $n$ (positive) parts of an integer of size at most $d$,
$$\alpha=(\alpha_1,\dots,\alpha_n), \ \ \sum_{i=1}^n\, \alpha_i\le d,$$ 
\item[(ii)] an unordered partition $\alpha'$ (with positive parts) of
$d-\sum_{i=1}^n\,\alpha_i$.
\end{enumerate}
The partition $\alpha'$ may be empty.
Let $\ell(\bar{\alpha})$ denote the length of $\bar{\alpha}$, 
$$\ell(\bar{\alpha})=\ell(\alpha)+\ell(\alpha')=n+\ell(\alpha').$$
Let $\Pi(d,n,k)$ be the set of partially ordered partitions of degree
$d$, order $n$, and length at least $d-k$. 
Our basic relations in the connected case will be indexed by $\Pi(d,n,k)$.

The set $\Pi(d,n,k)$ is stable for $k\ge d-n$: 
$$\Pi(d,n,k)=\Pi(d,n,d-n).$$ 
 Let $\Pi(d,n,\infty)$ denote the
stabilized set.

\subsubsection{The push-forward construction}
Let $g>0$ be the domain genus.
Let $\bar{\alpha}\in \Pi(d,n,k)$.
We will construct a relation 
$$T_{g, \bar{\alpha}}(\mu^1, \ldots, \mu^m \ | \ \gamma),$$ 
where:
\begin{enumerate}
\item[(i)] $\mu^1, \ldots, \mu^m$ are partitions of $d$,
\item[(ii)] $\gamma$ is a monomial $\prod_{j=1}^{n'}   \psi_{p_j}^{r_j}
             \cup \prod_{j=1}^m   \psi_{q_{j}}^{s_{j}}$. 
\end{enumerate}
The relation $T$ will be obtained from equivariant localization
        on 
\begin{equation}
\label{jjoo}
\overline{M}^\dagger_{g,n+\ell(\alpha'')+n'}(\mu^1,\ldots,\mu^m),
\end{equation}
where $\alpha''$ is the (possibly empty)
subpartition of $\alpha'$ consisting of parts of
size at least 2. We will use the abbreviated notation
$\overline{M}^\dagger$ for the moduli space of stable
relative maps (\ref{jjoo}).

The moduli space $\overline{M}^\dagger$ carries
$n+\ell(\alpha'')+n'$ markings of type $p$ indexed by the
following conventions:
\begin{enumerate}
\item[(i)] let $p_i$ denote the first $n$,
\item[(ii)] let $p_{i''}$ denote the middle $\ell(\alpha'')$,
\item[(iii)] and let $p_j$ denote the last $n'$.
\end{enumerate}
Let the middle $\ell(\alpha'')$ markings be placed in correspondence
with the parts of $\alpha''$.

Define the equivariant
relative Gromov-Witten class $\omega$ on the moduli space
$\overline{M}^\dagger$ as a product of four factors:
$$
        \omega =  \prod_{i=1}^n \psi_{p_i}^{\alpha_i-1} {\text{ev}}_{p_i}^*([\infty])
 \ \cdot 
      \prod_{i''=1}^{\ell(\alpha'')} \psi_{p_{i''}}^{\alpha''_{i''}-1} 
        {\text{ev}}_{p_{i''}}^*([\infty]) \ \cdot  \ \gamma \ \cdot $$
$$  ({\text{ev}}_{q_1}^*([0]))^{2+  \ell(\bar{\alpha})-d+k}  ,$$
where $\gamma$ is the monomial in the argument of $T$.    
The degree of $\omega$ is
          $$k + n +\ell(\alpha'') + \sum_{j=1}^{n'} r_j + 
\sum_{j=1}^m s_j + 2 .$$
Since the virtual dimension of 
$\overline{M}^\dagger$ is
          $$2g-2 + 2d + n+\ell(\alpha'') + n' + \sum_{i=1}^m (1+l(\mu_i)-d),$$
the dimension of $\omega \cap [\overline{M}^\dagger]^{vir}$ is
              $$  
2g-2 + 2d + n' + \sum_{i=1}^m (1+l(\mu_i)-d)
                                    -(k + \sum_{j=1}^{n'} r_j 
+ \sum_{j=1}^m s_j + 2) .$$

The moduli space $\overline{M}^\dagger$ carries $n+\ell(\alpha'')+n' + \sum_{i=1}^m
\ell(\mumu^i)$ total domain markings. Let
$$\rho'': \overline{M}^\dagger \rarr \overline{M}_{n+ n'+ \sum_{i=1}^m
\ell(\mu^i)}$$
be the stabilization map obtained from $\rho$ after forgetting the
markings corresponding to $\alpha''$.
Since the class $\text{ev}_{q_1}^*([0])$ occurs in $\omega$ with
total exponent at least 2, the non-equivariant limit of $\omega$ is 0.
Hence, in the non-equivariant limit, 
the push-forward $\rho''_*(\omega\cap [\overline{M}^\dagger]^{vir})$ vanishes:
\begin{equation}
\label{mainn}\rho''_* ( \omega \cap [\overline{M}^\dagger]^{vir}) = 0
\ \ \in A_*( \overline{M}_{g, n+ n'+ \sum_{i=1}^m
\ell(\mu^i)}).
\end{equation}
The left side of (\ref{mainn}) can be calculated by the localization
formula (\ref{ffkk}). As $\omega$ is a non-trivial equivariant class,
equation (\ref{mainn}) yields a non-trivial relation after
localizing and taking the non-equivariant limit,
 $$T_{g,{\bar {\alpha}}}(\mu^1, \ldots, \mu^m \ | \ \gamma)=
\sum_\Gamma \frac{m(\Gamma)}{|{\mathbf A}_\Gamma|} \ \rho''_{\Gamma*}\Big( \frac{\tau_\Gamma^*(\omega)} {e(N_\Gamma^{vir})}  \cap 
[\M_\Gamma]^{vir} \Big)=0,$$
where $\rho''_\Gamma= \rho'' \circ \tau_\Gamma.$

\subsubsection{The principal terms of $T$}

A localization graph $\Gamma$ corresponding
to a fixed locus of $\overline{M}^\dagger$ is of {\em type} $\bar{\beta}\in 
\Pi(d,n,k)
$
if the following properties are satisified:
\begin{enumerate}
\item[(i)] the vertex set $V_0$ consists of a single vertex $v_0$ of genus $g$,
\item[(ii)] the edge set $E$ is in bijective correspondence to the $\ell(\bar{\beta})$
parts of $\bar{\beta}$,
\item[(iii)] the vertex set $V_\infty$ consists of $\ell(\bar{\beta})$ vertices each
incident to a unique edge,
\item[(iv)] the $n+\ell(\alpha'')+ n'$ markings of type $p$
are distributed by the following rules:
\begin{enumerate}
\item[(a)] the first $n$ markings lie over $\infty$ with the $i^{th}$ marking
incident to the vertex corresponding to the  $i^{th}$ part of the
first partition $\beta$ of $\bar{\beta}$,
\item[(b)] the second $\ell(\alpha'')$ markings lie over $\infty$
on distinct vertices incident to the edges corresponding to the
parts of the second partition $\beta'$ of $\bar{\beta}$,
\item[(c)] the third $n'$ markings are all  incident to $v_0$,
\end{enumerate}
\item[(v)]
All refinements lie on side $0$.
\end{enumerate}
Let $\Gamma_{\bar{\beta}}$ denote the set of localization graphs
of type $\bar{\beta}$. Since $g>0$, the localization graphs
$\Gamma \in \Gamma_{\bar{\beta}}$ are {\em never} degenerate 
over 0.

The {\em principal} terms of the relation $T$ are indexed by partially
ordered partitions $\bar{\beta} \in \Pi(d,n,k)$.
The principal term of $T$ of  {\em type} $\bar{\beta}$ is:
$$T_{g,{\bar {\alpha}}}(\mu^1, \ldots, \mu^m \ | \ \gamma)[\bar{\beta}]=
\sum_{\Gamma\in \Gamma_{\bar{\beta}}} \frac{m(\Gamma)}{|{\mathbf A}_\Gamma|} \ \rho''_{\Gamma*}\Big( \frac{\tau_\Gamma^*(\omega)} {e(N_\Gamma^{vir})}  \cap 
[\M_\Gamma]^{vir} \Big).$$

We may compute the principal term of $T_{g,{\bar {\alpha}}}(\mu^1, \ldots, \mu^m \ | \ \gamma)$ 
of type $\bar{\beta}$
explicitly. However, we must define a function $S[\alpha''](\beta')$
 which arises naturally from the localization formula.

Let $\alpha''$ and $\beta'$ be two unordered partitions. Select
an ordering of the parts,
$$\alpha''=(\alpha''_1, \ldots, \alpha''_{\ell(\alpha'')}), \ \
\beta'=(\beta'_1, \ldots, \beta'_{\ell(\beta')}).$$
Define the integer sets  $\underline{0}=\emptyset$ and
$\underline{i}=\{1,\dots,i\}$ for each 
$i\ge1$.
Define the function
$S[\alpha''](\beta')$ by the following rules:
\begin{enumerate}
\item[(i)] if $\ell(\alpha'')>\ell(\beta')$, then $S[\alpha''](\beta')=0$,
\item[(ii)] if $\ell(\alpha'')\le \ell(\beta')$, then $S[\alpha''](\beta')$ 
is a sum over all injections
$$\iota:\underline{\ell(\alpha'')} \to\underline{\ell(\beta')},$$
$$
S[\alpha''](\beta')=\sum_{\iota:\underline{\ell(\alpha'')}
\to\underline{\ell(\beta')}}\,
\prod_{j=1}^{\ell(\alpha'')}\,\frac{1}{{\beta'_{\iota(j)}}^{\alpha''_j-1}} \cdot
\prod_{i\notin\im(\iota)}\,\frac{1}{\beta'_i}\,.
$$
\end{enumerate}
The function $S[\alpha''](\beta')$ depends only upon the {\em unordered}
partitions $\alpha''$ and $\beta'$.

Let $\overline{M}(\bar{\beta})$ denote the
moduli space of unparameterized maps,
$$\overline{M}(\bar{\beta})=
\overline{M}_{g,n'}(\bar{\beta},\mu^1, \ldots, \mu^m),$$
with markings of type $q$ indexed by $0,1,\ldots,m$.
Define the standard Gromov-Witten class
 $\omega_{\bar{\beta}}$ on the moduli space $\overline{M}(\bar{\beta})$
by:
$$\omega_{\bar{\beta}} =  \gamma \ \cdot \
                      \psi_{q_0}^{\ell(\bar{\beta})-d+k}.$$
The degree of $\omega_{\bar{\beta}}$ is
$$  \sum_{j=1}^{n'} r_j + \sum_{j=1}^m s_j +\ell(\bar{\beta}) -d +k  .$$
Since the dimension of $\overline{M}(\bar{\beta})$ is
$$2g-2 + 2d + n' + (1+\ell(\bar{\beta})-d) +\sum_{j=1}^m 
(1+\ell(\mu^j)-d) -3,$$
the dimension of
$\omega_{\bar{\beta}} \cap [\overline{M}(\bar{\beta})]^{vir}$ is
$$
2g-2 + 2d + n' + \sum_{i=1}^m (1+l(\mu_i)-d)
                                    -(k + \sum_{j=1}^{n'} r_j 
+ \sum_{j=1}^m s_j + 2),$$
equal to the dimension of $\omega \cap [\overline{M}^\dagger]^{vir}$.

A direct application of the localization formula (\ref{ffkk}) yields
the following result.

\begin{lm} The principal term of 
$T_{g,{\bar {\alpha}}}(\mu^1, \ldots, \mu^m \ | \ \gamma)$
of type $\bar{\beta}$ is
$$\Big( \prod_{i=1}^n \frac{1}{\beta_i^{\alpha_i-1}} \ S[\alpha''](\beta')
\  (-1)^{n+\ell(\beta')} \ \eta(\bar{\beta})\Big) \ \cdot \
\rho''_*(\omega_{\bar{\beta}} \cap [\overline{M}(\bar{\beta})]^{vir}),$$ 
where $\eta(\bar{\beta})$ is the non-vanishing factor
$$\frac{1}{|\text{\em Aut}(\beta')|}
\prod_{i=1}^n \frac{1}{(-1)^{\beta_i} \frac{\beta_i!}{\beta_i^{\beta_i}}}
\prod_{i=1}^{\ell(\beta')} \frac{1}{(-1)^{\beta'_i} \frac{\beta'_i!}{{\beta'_i}^{\beta'_i}}}.
$$
\end{lm}

\subsection{Basic relations II: The disconnected case}

\label{dizzz}

\subsubsection{The set $\Pi(\DD,\NN,k)$}
\label{pidnk00}  
Let $\DD=(d_1,\ldots,d_c)$ be an ordered degree partition, and let
 $$d=\sum_{i=1}^c d_i.$$
Let
$\NN=(n_1,\ldots,n_c)$ be ordered set partition of $\underline{n}$
for which
$$d_i\ge |n_i| >0,$$
for all $i$.
Let $k\ge0$ be an integer.
A partially ordered partition,
 $$\bar{\alpha}=((\alpha[1], {\alpha'[1]}),\ldots,
(\alpha[c],{\alpha'[c]})),$$ of degree $\DD$ and order $\NN$
consists of the following data:
\begin{enumerate}
\item[(i)] ordered
partitions with $|n_i|$ (positive) parts of integers of size at most $d_i$,
$$\alpha[i]=(\alpha[i]_1,\dots,\alpha[i]_{|n_i|}), 
\ \ \sum_{j=1}^n\, \alpha[i]_j\le d_i,$$ 
\item[(ii)] unordered partitions ${\alpha'[i]}$ (with positive parts) of
$d_i-\sum_{j=1}^n\,\alpha[i]_j$.
\end{enumerate}
The partitions ${\alpha'[i]}$ may be empty.
Let $\ell(\bar{\alpha})$ denote the length of $\bar{\alpha}$, 
$$\ell(\bar{\alpha})=\sum_{i=1}^c (\ell(\alpha[i])+\ell({\alpha'[i]}) )
=n+\sum_{i=1}^c \ell({\alpha'[i]}).$$
Let $\Pi(\DD,\NN,k)$ be the set of partially ordered partitions of degree
$\DD$, order $\NN$, and length at least $d-k$. 
Our basic relations in the disconnected case will be indexed by $\Pi(\DD,\NN,k)$.

The set $\Pi(\DD,\NN,k)$ is stable for $k\ge d-n$: 
$$\Pi(\DD,\NN,k)=\Pi(\DD,\NN,d-n).$$ 
 Let $\Pi(\DD,\NN,\infty)$ denote the
stabilized set.

\subsubsection{The push-forward construction}
Let $\GG=(g_1, \ldots,g_c)$ be an ordered set of genera where not all $g_i=0$.
Let $\bar{\alpha}\in \Pi(\DD,\NN,k)$.
We will construct a relation 
$$T_{\GG, \bar{\alpha}}(\mumu^1, \ldots, \mumu^m \ | \ \gamma),$$ 
where:
\begin{enumerate}
\item[(i)] $\mu^1, \ldots, \mu^m$ are partitions of type $\DD$,
\item[(ii)] $\gamma$ is a monomial $\prod_{j=1}^{n'}   \psi_{p_j}^{r_j}
             \cup \prod_{j=1}^m   \psi_{q_{j}}^{s_{j}}$. 
\end{enumerate}
The relation $T$ will be obtained from equivariant localization
        on 
\begin{equation}
\label{jjjoo}
\overline{M}^\dagger_{\GG,\NN+\ell(\alpha'')+\NN'}(\mumu^1,\ldots,\mumu^m),
\end{equation}
where $\alpha''$ is the (possibly empty)
subpartition of $({\alpha'[1]}, \ldots, {\alpha'[c]})$ consisting of parts of
size at least 2. The last $n'$ markings of type $p$ are distributed on the
domain components by $\NN'$.
As before, let 
$\overline{M}^\dagger$ denote the moduli space of stable
relative maps (\ref{jjjoo}).

Define the equivariant relative
Gromov-Witten class $\omega$ on the moduli space
$\overline{M}^\dagger$ as a product of four factors:
\begin{equation}
\label{yepp}        \omega  =  \prod_{i=1}^c \Big( \prod_{j=1}^{|n_i|} \psi_{p_j}^{\alpha[i]_j-1} {\text{ev}}_{p_j}^*([\infty])
 \ \cdot 
      \prod_{j''=1}^{\ell(\alpha''[i])} \psi_{p_{j''}}^{\alpha''[i]_{j''}-1} 
        {\text{ev}}_{p_{j''}}^*([\infty]) \Big) \ \cdot  \ \gamma \ \cdot 
\end{equation}
$$({\text{ev}}_{q_1}^*([0]))^{ 2+\ell(\bar{\alpha})-d+k  },$$
where $\gamma$ is the monomial in the argument of $T$.   

The dimension calculus proceeds exactly as in the connected case (replacing
the connected genus $g$ by the arithmetic genus $\sum_{i=1}^c g_i -c+1$
in the disconnected case).

    The moduli space $\overline{M}^\dagger$ carries $n+\ell(\alpha'')+n' + \sum_{i=1}^m
\ell(\mu^i)$ total domain markings. Let
$$\rho'': \overline{M}^\dagger \rarr \prod_{i=1}^c
\overline{M}_{g_i,|n_i|+ |n'_i|+ \sum_{j=1}^m
\ell(\mu^j[i])}$$
be the stabilization map obtained from $\rho$ after forgetting the
markings corresponding to $\alpha''$.
As before, the push-forward vanishes in the non-equivariant limit,
\begin{equation}
\label{mainnn}\rho''_* ( \omega \cap [\overline{M}^\dagger]^{vir}) = 0.
\end{equation}
The left side of (\ref{mainnn}) can be calculated by the localization
formula (\ref{ffkk}) to yield the relation
 $$T_{\GG,{\bar {\alpha}}}(\mumu^1, \ldots, \mumu^m \ | \ \gamma)=
\sum_\Gamma \frac{m(\Gamma)}{|{\mathbf A}_\Gamma|} \ \rho''_{\Gamma*}\Big( \frac{\tau_\Gamma^*(\omega)} {e(N_\Gamma^{vir})}  \cap 
[\M_\Gamma]^{vir} \Big)=0.$$

\subsubsection{The principal terms of $T$}

A localization graph $\Gamma$ corresponding
to a fixed locus of $\overline{M}^\dagger$ is of {\em type} $\bar{\beta}\in 
\Pi(\DD,\NN,k)
$
if the following properties are satisified:
\begin{enumerate}
\item[(i)] the vertex set $V_0$ consists of vertices $v_{1,0},\ldots
v_{c,0}$ with genus assignments $g_1, \ldots, g_c$ respectively,
\item[(ii)] the edge set $E$ is in bijective correspondence to the $\ell(\bar{\beta})$
parts of $\bar{\beta}$,
\item[(iii)] the vertex $v_{i,0}$ is incident to the edges corresponding to
the parts of $(\beta[i], {\beta'[i]})$,
\item[(iv)] the vertex set $V_\infty$ consists of $\ell(\bar{\beta})$ vertices each
incident to a unique edge,
\item[(v)] the $n+\ell(\alpha'')+ n'$ markings of type $p$
are distributed by the following rules:
\begin{enumerate}
\item[(a)] the markings $n_i$ lie over $\infty$ with the $j^{th}$ marking
incident to the vertex corresponding to the  $j^{th}$ part of 
$\beta[i]$,
\item[(b)] the middle $\ell({\alpha''[i]})$ markings lie over $\infty$
on distinct vertices incident to the edges corresponding to the
parts of ${\beta'[i]}$,
\item[(c)] the markings $n_i'$  are   incident to $v_{i,0}$,
\end{enumerate}
\item[(vi)]
All refinements lie on side $0$.
\end{enumerate}
Let $\Gamma_{\bar{\beta}}$ denote the set of localization graphs
of type $\bar{\beta}$.
Since not all $g_i=0$, the localization graphs
$\Gamma \in \Gamma_{\bar{\beta}}$ are {\em never} degenerate 
over 0.

The {\em principal} terms of the relation $T$ are indexed by partially
ordered partitions $\bar{\beta} \in \Pi(\DD,\NN,k)$.
The principal term of $T$ of  {\em type} $\bar{\beta}$ is:
$$T_{\GG,{\bar {\alpha}}}(\mumu^1, \ldots, \mumu^m \ | \ \gamma)[\bar{\beta}]=
\sum_{\Gamma\in \Gamma_{\bar{\beta}}} \frac{m(\Gamma)}{|{\mathbf A}_\Gamma|} \ \rho''_{\Gamma*}\Big( \frac{\tau_\Gamma^*(\omega)} {e(N_\Gamma^{vir})}  \cap 
[\M_\Gamma]^{vir} \Big).$$

The principal term of $T_{g,{\bar {\alpha}}}(\mu^1, \ldots, \mu^m \ | \ \gamma)$ 
of type $\bar{\beta}$
may also be explicitly computed
from the localization formula.

Let $\overline{M}(\bar{\beta})$ denote the
moduli space of unparameterized maps,
$$\overline{M}(\bar{\beta})=
\overline{M}_{\GG,\NN'}(\bar{\beta},\mumu^1, \ldots, \mumu^m),$$
with markings of type $q$ indexed by $0,1,\ldots,m$.
Define the standard Gromov-Witten class
 $\omega_{\bar{\beta}}$ on the moduli space $\overline{M}(\bar{\beta})$
by:
$$\omega_{\bar{\beta}} =  \gamma \ \cdot \
                      \psi_{q_0}^{\ell(\bar{\beta})-d+k}.$$

A direct application of the localization formula (\ref{ffkk}) yields
the following result.

\begin{lm}
\label{coey}
 The principal term of 
$T_{\GG,{\bar {\alpha}}}(\mumu^1, \ldots, \mumu^m \ | \ \gamma)$
of type $\bar{\beta}$ is
$$\prod_{i=1}^c \Big( 
\prod_{j=1}^{|n_i|} \frac{1}{(\beta[i]_j)^{\alpha[i]_j-1}} \ 
S[{\alpha''[i]}]({\beta'[i]})
\  (-1)^{|n_i|+\ell({\beta'[i]})} \ \eta((\beta[i],{\beta'[i]})) \Big) $$
$$ \cdot \
\rho''_*(\omega_{\bar{\beta}} \cap [\overline{M}(\bar{\beta})]^{vir}).$$
\end{lm}

\subsection{The matrix ${\mathbf M}$}
Let ${\mathbf M}(\DD,\NN,k)$ be the  matrix with rows and columns
indexed by the set $\Pi(\DD,\NN,k)$ defined by the
prefactors of the principal terms of the basic relations. For 
$$\bar{\alpha}, \bar{\beta} \in \Pi(\DD,\NN,k),$$
let the element of ${\mathbf M}(\DD,\NN,k)$ in position $(\bar{\alpha}, 
\bar{\beta})$
be:
$$\prod_{i=1}^c \Big( 
\prod_{j=1}^{|n_i|} \frac{1}{(\beta[i]_j)^{\alpha[i]_j-1}} \ 
S[{\alpha''[i]}]({\beta'[i]})
\  (-1)^{|n_i|+\ell({\beta'[i]})} \ \eta((\beta[i],{\beta'[i]})) \Big),$$
the prefactor of the
principal term of 
$T_{\GG,{\bar {\alpha}}}(\mumu^1, \ldots, \mumu^m \ | \ \gamma)$
of type ${\bar{\beta}}$.

\begin{lm}
\label{mmm}
${\mathbf M}(\DD,\NN,k)$ is invertible.
\end{lm}

\noindent

Lemma \ref{mmm} will be proven in Section \ref{matt}. The nonsingularity
of ${\mathbf M}(\DD,\NN,k)$ plays a fundamental role in our proof
of Theorem 2 part (ii).

\subsection{Proof of Theorem 2 part (ii)}
\subsubsection{Overview}
Let $\overline{M}_{\GG,\NN'}(\mumu^1, \ldots, \mumu^m)$ be a
moduli space of stable relative maps with $c$ domain components,
degree partition $\DD$, and total degree $d$.
Let $\omega$ be a relative Gromov-Witten class on 
 $\overline{M}_{\GG,\NN'}(\mumu^1, \ldots, \mumu^m)$.
Let $$\rho_*(\omega\cap [\overline{M}]^{vir})$$
denote the push-forward
$$\rho_*\Big( \omega \cap  [\overline{M}_{\GG,\NN'}(\mumu^1, \ldots, 
\mumu^m)]^{vir}\Big).$$
We must prove $\rho_*(\omega\cap [\overline{M}]^{vir})$ 
is a tautological class on the
target of $\rho$.

The proof consists of a multilayer induction using the basic
relations $T$. We will
prove the result for each total degree $d$ separately.
The proof will follow the notation of Section \ref{dizzz}.

\subsubsection{Genus induction}
Since all algebraic cycle classes on products of moduli spaces of
genus 0 pointed curves are
tautological [K], the class $$\rho_*(\omega\cap [\overline{M}]^{vir})$$
 is certainly
tautological if all the domain genera $g_i$ are 0.

We will proceed by induction on the arithmetic genus,
$$g(\GG)= \sum_{i=1}^c g_i -c +1,$$
of the domain curve.
The smallest value for the arithmetic genus is
$-d+1$ which occurs when the number of domain components equals $d$ and
$$\GG=(0,\ldots,0), \ \  \DD=(1, \ldots,1).$$
The base case of the genus induction holds as all domain genera
in the base case are 0.

\subsubsection{Induction on the markings}
We now induct upon the 
number of markings $m$ of type $q$. We first establish the
base case $m=2$. In fact, all aspects of the proof of the full result are
manifest in the base case.
 
We will also induct on the number $n'$ of markings of type $p$. The
induction on type $p$ markings is straightforward and plays a very minor role.

Let $\GG$ be a genus distribution with not all genera 0, and
let $\DD$ be a degree distribution.
Let a partition $\mu^1$ of type $\DD$ and a relative Gromov-Witten class,
$$\gamma=  
\prod_{j=1}^{n'}   \psi_{p_j}^{r_j}
             \cup \psi_{q_{1}}^{s_{1}},$$ 
be given.

Let $\NN$ have cardinalities $|\NN|=(1, \ldots,1)$. 
We will first consider all basic relations,
\begin{equation}
\label{rrr}
T_{\GG,{\bar {\alpha}}}(\mumu^1\ | \ \gamma),
\end{equation}
where $\bar{\alpha}\in \Pi(\DD,\NN,k)$.

Let us now analyze the terms of the relations (\ref{rrr}).
Let $\Gamma$ be a localization graph. Since $\text{ev}_{p_j}^*([\infty])$
occurs in  (\ref{yepp}) for all the first $n+\ell(\alpha'')$ markings
of type $p$, only graphs $\Gamma$ for which these markings lie over $\infty$
contribute non-trivially to the basic relations.
Similiarly,
since $\text{ev}_{q_1}^*([0])$
occurs in (\ref{yepp}), only graphs $\Gamma$ for which
$R^1$ lies on side 0 contribute.

As there are no refinements on side $\infty$, we are in Case II
of the localization formula. We find,
$$\overline{M}_\Gamma = \overline{M}_0 \times \overline{M}_\infty,$$
where
$$\overline{M}_0= \overline{M}_{\GG_0, \NN'_0}(R^1, R_\delta)$$
and
$\overline{M}_\infty$ is simply a product of
moduli spaces of pointed curves.
If
\begin{equation}
\label{pwp} 
g(\GG_0) < g(\GG) \ \  or \ \ n'_0 < n',
\end{equation} then the
localization terms over $0$ of $\Gamma$ push-forward under $\rho''$ to
tautological classes by the induction hypotheses on the genus and
the marking number $n'$. 
The localization terms over $\infty$ certainly push-forward to
tautological classes.

The graphs $\Gamma$ satisfying $g(\GG_0)=g(\GG)$ have only genus 0,
edge valence 1 vertices over $\infty$. If $n'_0=n'$, only the
first $n+\ell(\alpha'')$ markings of type $p$ lie over $\infty$. 
Of these, the first $n$ markings are distributed on the different connected
components of $\Gamma$ since $|\NN|=(1,\ldots,1)$.
Since all the parts of
$\alpha''$ are at least two,  the second $\ell(\alpha'')$ markings of
type $p$ trivialize the contribution of $\Gamma$ if {\em any}
vertex of marking valence greater than 1 occurs over $\infty$.
An ordered partition $\bar{\beta}$ is obtained from $\Gamma$ from
the edge degrees. By dimension considerations, the contribution of
$\Gamma$ is trivial unless $\bar{\beta} \in \Pi(\DD,\NN,k)$.

We have proven, for contributing graphs $\Gamma$, {\em either} condition
(\ref{pwp}) holds {\em or} 
$$\Gamma \in \Gamma_{\bar{\beta}}.$$ 
Therefore, the non-principal terms of the basic relations (\ref{rrr})
all push-forward via $\rho''$ to tautological classes.

The principal terms of the basic relations (\ref{rrr}) taken for
all $\bar{\alpha} \in \Pi(\DD,\NN,k)$ determine a nonsingular matrix
of prefactors 
by Lemmas \ref{coey} and \ref{mmm}.  We conclude, for all 
$\bar{\beta} \in \Pi(\DD,\NN,k)$,
\begin{equation}
\label{ppqq}
\rho''_*(\gamma \cdot \psi_{q_0}^{\ell(\bar{\beta})-d+k} \cap [\overline{M}_{\GG,\NN'}(\bar{\beta}, \mumu^1)]^{vir})
\end{equation}
is tautological in the target moduli of curves.
By considering the classes (\ref{ppqq}) for all $k\geq 0$,
we find,
\begin{equation}
\label{ppqqqq}
\rho''_*(\omega
 \cap [\overline{M}_{\GG,\NN'}(\bar{\beta}, \mumu^1)]^{vir})
\end{equation}
is tautological on the target for all relative Gromov-Witten
classes $\omega$ and all $\bar{\beta} \in \Pi(\DD,\NN,\infty)$. 
We have exhausted the basic relations
(\ref{rrr}) for the element $|\NN|=(1,\ldots,1)$.

If $\DD=(1,\ldots,1)$, we have completed the $m=2$ induction base 
since $\rho =\rho''$.
However, if $\DD$ has parts of size greater
than $1$, then $\rho \neq \rho''$ for general $\bar{\beta}$. 
To proceed, consider an augmentation
$$|\NN|=(1,\ldots,1,2,1,\ldots,1),$$
subject to $d_i\geq |n_i| >0$. 
We will exhaust the basic relations for the augmented $\NN$.
We will repeat the cycle of augmentation and exhaustion until
all $\NN$ satisfying $d_i\geq |n_i| >0$ have been considered. 

We repeat the entire argument for the augmented $\NN$.
The argument is identical except for one important difference.
For graphs $\Gamma$ satisfying 
$$g(\GG_0)=g(\GG), \ \ n'_0=n',$$
the first $n$ markings of type $p$ may {\em not} be incident to
distinct vertices over $\infty$. However, the graphs which
contain common incidences yield tautological push-forwards
by the analysis for {\em lower} $\NN$. 
After each cycle of the analysis, we conclude,
\begin{equation}
\label{ppqqq}
\rho''_*(\omega
 \cap [\overline{M}_{\GG,\NN'}(\bar{\beta}, \mumu^1)]^{vir})
\end{equation}
is tautological in the target for all relative Gromov-Witten
classes $\omega$ and $\bar{\beta} \in \Pi(\DD,\NN,\infty)$.
The $m=2$ base case is proven.

\subsubsection{The full induction}
The full induction exactly follows the argument for the
$m=2$ base case.
Let $m\geq 3$.

Let $\GG$ be a genus distribution with not all genera 0, and
let $\DD$ be a degree distribution.
Let partitions $\mumu^1, \ldots, \mumu^{m-1}$ of type 
$\DD$ and a relative Gromov-Witten class,
$$\gamma=  
\prod_{j=1}^{n'}   \psi_{p_j}^{r_j}
             \cup \prod_{j=1}^{m-1} \psi_{q_{j}}^{s_{j}},$$
be given.

We consider the basic relations,
\begin{equation}
\label{rrrr}
T_{\GG,{\bar {\alpha}}}(\mumu^1, \ldots, \mumu^{m-1}\ | \ \gamma),
\end{equation}
where $\bar{\alpha}\in \Pi(\DD,\NN,k)$.

Let $\Gamma$ be a localization graph. 
As before, only graphs $\Gamma$  for which the first $n+\ell(\alpha'')$ markings
of type $p$ lie over $\infty$
contribute non-trivially to the basic relations.
Since $\text{ev}_{q_1}^*([0])$
occurs in (\ref{yepp}), only graphs $\Gamma$ for which
$R^1$ lies on side 0 contribute.

If there are refinements over $\infty$, then we are in Case I
of the localization formula.
We find,
$$\overline{M}_\Gamma = \overline{M}_0 \times \overline{M}_\infty,$$
where
$$\overline{M}_0= \overline{M}_{\GG_0, \NN'_0}({\mathbf R}_0, R_\delta), \ \
\overline{M}_\infty= \overline{M}_{\GG_0, \NN'_0}(R_\delta, {\mathbf R}_\infty). 
$$
Since the graph $\Gamma$ has at most ${m-2}$ refinements on side $\infty$,
both $\overline{M}_0$ and $\overline{M}_\infty$ have strictly {\em less}
than $m$ markings of type $q$.
The
localization terms over $0$ and $\infty$ of $\Gamma$ push-forward under $\rho''$ to
tautological classes by the induction hypotheses on $m$.

If there are no refinements over $\infty$,
 we are in Case II
of the localization formula. 
We then exactly follow the augmentation and
exhaustion argument of the base case to conclude:
\begin{equation}
\label{ppqqqqq}
\rho''_*(\omega
 \cap [\overline{M}_{\GG,\NN'}(\bar{\beta}, \mumu^1, \ldots, \mumu^{m-1})]^{vir})
\end{equation}
is tautological on the target for all relative Gromov-Witten
classes $\omega$ and all $\bar{\beta} \in \Pi(\DD,\NN,\infty)$
for all $\NN$ subject to
 $d_i\geq |n_i| >0$. 
The induction step is therefore established and the proof
of Theorem 2 part(ii) is complete.
\epf

\subsubsection{Effectivity}
The argument of Theorem 2 part (ii) yields an effective
procedure for calculating the class of the
push-forward,
\begin{equation}
\label{wxx}
\rho_*(\omega \cap [\overline{M}]^{vir}),
\end{equation}
in the tautological ring of the target.

To establish the procedure, we must only check the push-forward
(\ref{wxx}) can be effectively determined if
\begin{equation}
\label{zzw}
\GG=(0, \ldots, 0),
\end{equation}
the base case of the genus induction.
In fact, the argument of Theorem 2 part (ii) applies with almost
no modification
to the base case (\ref{zzw}).
A minor 
issue is the possibility of degenerate data over the fixed point $0$ (for
which the dimension formula is not valid).
We leave the details to the interested reader. The outcome
is a procedure for calculating (\ref{wxx}) in the base case (\ref{zzw}).

Theorem 2 part (i) is then also effective. The push-forward,
$$\rho_*(\omega \cap [\overline{M}^\dagger]^{vir}),$$
is expressible by localization in terms of the
push-forwards of relative Gromov-Witten classes on
 unparameterized spaces.

\section{Nonsingularity}
\label{matt}
\subsection{Overview} 
Let $\mata(\DD,\NN,k)$ be the matrix with rows and columns indexed
by $\Pi(\DD,\NN,k)$ with the following coefficients.
For $$\bar{p}, \bar{q} \in \Pi(\DD,\NN,k),$$
let the element of $\mata(\DD,\NN,k)$ in position $(\bar{p}, 
\bar{q})$
be:
$$\prod_{i=1}^c \Big( 
\prod_{j=1}^{|n_i|} \frac{1}{(p[i]_j)^{q[i]_j-1}} \ 
S[{q''[i]}]({p'[i]}) \Big).$$
The matrix $\mata(\DD,\NN,k)$ is obtained from ${\mathbf M}(\DD,\NN,k)$
by transposition followed by a rescaling of the rows and columns.
We will prove Lemma \ref{mmm} by establishing the nonsingularity
of $\mata(\DD,\NN,k)$.

\subsection{Nonsingularity for $\Pi(d,n,k)$}
\subsubsection{Overview}
We will start by proving the nonsingularity of $\mata$ in
case the index set is $\Pi(d,n,k)$.

\begin{lm}
\label{lmaa}
The matrix $\mata(d,n,k)$ is invertible.
\end{lm}

The nonsingularity of $\mata(\DD,\NN,k)$ will be a direct
consequence of the proof of Lemma \ref{lmaa}.

\subsection{The function $T[q''](p'')$}

In order to prove Lemma \ref{lmaa}, we will require an auxiliary matrix.
Let $p''$ and $q''$ be partitions with parts of size at least 2. 
Define the function $T[q''](p'')$ by the following rules:
\begin{enumerate}
\item[(i)] If $\ell(q'')>\ell(p'')$, then $T[q''](p'')=0$.
\item[(ii)] If $\ell(q'')\le \ell(p'')$, 
then $T[q''](p'')$ is a sum over all injections
$$\theta:\underline{\ell(q'')}\to\underline{\ell(p'')}.$$ 
For such an injection $\theta$,
let 
$$v(\theta)=|\{j\in\underline{\ell(q'')}:\quad q''_j=p''_{\theta(j)}=2\}|.
$$
Then,
$$
T[q''](p'')=\sum_{\theta:\underline{\ell(q'')}\to\underline{\ell(p'')}}\,
2^{v(\theta)}\,\prod_{j=1}^{\ell(q'')}\,\binom{p''_{\theta(j)}-1}{q''_j-1}
(-1)^{q''_j-1}\,{(q''_j)}^{p''_{\theta(j)}-2}\,.
$$
\end{enumerate}

Let $\bar{p}, \bar{q} \in \Pi(d,n,k)$. 
Let $p''$ and $q''$ be the subpartitions of 
$p'$ and $q'$ respectively consisting of the parts of size at least $2$.
Let
$p'''$ be the subpartition of $p''$ consisting of parts at least $3$.

Let 
$\matb(d,n,k)$ be the matrix  with  rows and columns
indexed by the set $\Pi(d,n,k)$ with coefficient  
\begin{equation*}
\frac{1}{|\Aut(p'')|\,|\Aut(q'')|}\prod_{h}\,p_h\,\prod_i\,p'''_i\,
\prod_{j=1}^n\Big(\binom{p_j-1}{q_j-1}(-1)^{q_j-1}q_j^{p_j-2}\Big)
\cdot T[q''](p'')\,,
\end{equation*}
in position $(\bar{p}, \bar{q})$.

\subsection{An order on $\Pi(d,n,k)$}

We order the set
$\Pi(d,n,k)$ in the following manner. Let $\bar{p}$ and $\bar{q}$ be
distinct
elements of $\Pi(d,n,k)$. 
We assume the parts of the subpartitions
$p''$ and $q''$ are arranged in increasing order.
Then, $\bar{p}$ precedes $\bar{q}$ if and only if
\begin{enumerate}
\item[$\bullet$] $\ell(p'')<\ell(q'')$, or
\item[$\bullet$] $\ell(p'')=\ell(q'')$ and $p''$ precedes $q''$ in the lexicographic
order, or
\item[$\bullet$] $p''=q''$ and $\ell(p')>\ell(q')$, or
\item[$\bullet$] $p'=q'$ and $p$ precedes $q$ in the lexicographic
order.
\end{enumerate}
We write
$\bar{p} < \bar{q}$
if $\bar{p}$ precedes $\bar{q}$ in the
above order on $\Pi(d,n,k)$.

\subsection{Proof of Lemma \ref{lmaa}}
\label{plmaa}
Define $\matc(d,n,k)$ by 
$$\matc(d,n,k)=\matb(d,n,k)  \cdot  \mata(d,n,k).$$
Lemma \ref{lmaa} is an immediate consequence of the following result.
\begin{lm}
\label{lmbb}
The matrix $\matc(d,n,k)$ is upper triangular
with $1$'s along the diagonal. 
\end{lm}
\bpf
Recall the well-known identities
$$\sum_{k=0}^n\binom{n}{k}(-1)^kk^a=0,\qquad0\le a\le n-1,$$
and
$$\sum_{k=0}^n\binom{n}{k}(-1)^k(k+1)^{-1}=\frac{1}{n+1}\,,\qquad n\ge0.$$
To prove these identities, consider the function 
$$f(t)=\sum_{k=0}^n\binom{n}{k}(-1)^kt^k=(1-t)^n.$$
The first identity is obtained
by the evaluation of $(t\frac{d}{dt})^a f$ at $t=1$, and the second
is obtained by computing $\int_0^1f(t)dt$. 

Let $\pti$, $\rti$, $\pbar$, and $\rbar$ be integers satisfying
$\pti\ge \rti\ge1$ and $\pbar\ge\rbar\ge2$.
Further, let $\delta(x,y,z)=1$ when $x=y=z$, and $0$ otherwise.
We will need four closely related sums:
\begin{eqnarray*}
(\alpha):\quad&&\sum_{\qti=1}^{\pti}\binom{\pti-1}{\qti-1}
(-1)^{\qti-1}\qti^{\pti-1-\rti}
=\begin{cases} 0, & \rti<\pti, \\ \frac{1}{\pti}\,,& \rti=\pti;\end{cases} \\
(\beta):\quad&&\sum_{\qbar=1}^{\pbar}\binom{\pbar-1}{\qbar-1}
(-1)^{\qbar-1}\qbar^{\pbar-3}2^{\delta(2,\pbar,\qbar)}
=\begin{cases} 0, &\pbar\ge3,\\ 0,&\pbar=2; \end{cases}\\
(\beta'):\quad&&\sum_{\qbar=1}^{\pbar}\binom{\pbar-1}{\qbar-1}
(-1)^{\qbar-1}\qbar^{\pbar-2}2^{\delta(2,\pbar,\qbar)}
=\begin{cases} 0, &\pbar\ge3,\\ -1,&\pbar=2; \end{cases}\\
(\gamma):\quad&&\sum_{\qbar=1}^{\pbar}\binom{\pbar-1}{\qbar-1}
(-1)^{\qbar-1}\qbar^{\pbar-1-\rbar}2^{\delta(2,\pbar,\qbar)}
=\begin{cases}  0, &\rbar<\pbar,\\ 0,& \rbar=\pbar=2,\\ \frac{1}{\pbar},&
\rbar=\pbar\ge3. \end{cases}
\end{eqnarray*}

Let $\bar{p}$ and $\bar{r}$ be elements of $\Pi(d,n,k)$. The matrix element
$\matc(\bar{p},\bar{r})$ is given by
\begin{equation}
\label{mss}
\matc(\bar{p},\bar{r})=\sum_{\bar{q}\in \Pi(d,n,k)} \,
\matb(\bar{p},\bar{q})\mata(\bar{q},\bar{r}) 
\end{equation}
$$=\frac{\prod\,p_h\,\prod\,p'''_i}{|\Aut(p'')|}\,\sum_{\bar{q}}
\prod_{j=1}^n\,\Big(\binom{p_j-1}{q_j-1}(-1)^{q_j-1}q_j^{p_j-1-r_j}\Big)
\cdot \frac{T[q''](p'')}{|\Aut(q'')|}\cdot S[r''](q'). $$
We must show:
\begin{enumerate}
\item[(i)]
 $\matc(\bar{p}, \bar{r})=0$ if $\bar{r}$ precedes
$\bar{p}$ in the order on $\Pi(d,n,k)$,
\item[(ii)]
 $\matc(\bar{p},\bar{p})=1$ for all $\bar{p}$.
\end{enumerate}

Let $\bar{r}\le \bar{p}$ with respect to the order  on $\Pi(d,n,k)$.
The
prefactor $\frac{\prod\,p_h\,\prod\,p'''_i}{|\Aut(p'')|}$ plays
a role only in the normalization of the diagonal elements (ii).
We will first concentrate on 
the main sum  in  formula (\ref{mss})
for $\matc(\bar{p},\bar{r})$: a triple sum over partially
ordered partitions $\bar{q}$, injections 
$$\theta:\underline{\ell(q'')}\to \underline{\ell(p'')}$$ via $T[q''](p'')$,
and injections 
$$\iota:\underline{\ell(r'')}\to \underline{\ell(q')}$$ via $S[r''](q')$.

The key step in the proof is to recognize that 
the main sum in (\ref{mss}) can be written
in terms of products of the sums $(\alpha)$, $(\beta)$, $(\beta')$,
and $(\gamma)$ discussed above. The products have $n+|p''|$ factors.
We will define and analyze the products in stages.

We first consider the case $p''=\emptyset$. We will write the
sum (\ref{mss}) as a single product with $n$  factors of 
type $(\alpha)$.
The $n$ factors are obtained by replacing
$(\pti,\qti,\rti)$ in the formula for
$(\alpha)$
by $(p_i,q_i,r_i)$, for $1\le i\le n$. As $q_i$ varies in the
range,
$$1 \leq q_i \leq p_i,$$
all possible choices
of the ordered partition, 
$$q=(q_1,\dots,q_n),$$ in the partially ordered
partition $\bar{q}=(q,q')$ indexing (\ref{mss}) occur, as
$\bar{q}$ does not contribute to (\ref{mss}) if there exists
an $i$ for which $q_i>p_i$.

Since $p''=\emptyset$, we find for $\bar{q}$ to contribute to (\ref{mss}),
the equality
$q''=\emptyset$ must hold.
By the
ordering $\bar{r} \leq \bar{p}$, the equality 
$r''=\emptyset$ must hold.
The sum over $\bar{q}$, $\theta$, and $\iota$
is then just a sum over $\bar{q}$. The partially ordered
partition $\bar{q}$ is given by the choice of
$q=(q_1,\dots,q_n)$, for $q'$ then consists of $d-\sum q_i$ parts equal
to $1$.  We conclude
the main sum of (\ref{mss}) equals the product of the
$n$ sums of type $(\alpha)$.
 
The product of the $n$ sums 
vanishes if there exists an $i$ for which
$r_i<p_i$. Since $\bar{r}\le \bar{p}$ and $r''=p''$, we find $|r'|\ge|p'|$.
Hence, $\sum_{i=1}^n r_i\le
\sum_{i=1}^n p_i$. For a nonvanishing product, $\bar{r}=\bar{p}$ is then 
necessary. In the diagonal case, the product equals $\prod_{i=1}^n (1/p_i)$ 
which is exactly cancelled by the prefactor
to yield the matrix element $1$. We have completed the
analysis of the case $p''=\emptyset$.

We now consider general $p''$.
Let $m=\ell(p'')$. We will analyze associated products with
$n$ factors
of type $(\alpha)$ and $m$ additional factors, each
of type
$(\beta)$, $(\beta')$, or $\gamma$.
The first $n$ factors are obtained as before by replacing
$(\pti,\qti,\rti)$ in the formula for
$(\alpha)$
by $(p_i,q_i,r_i)$, for $1\le i\le n$. 
Let the $m$ upper bounds of the last sums be denoted by
$\pbar_1,\dots,\pbar_m$, and let the $m$ index variables 
for the sums be
$\qbar_1,\dots,\qbar_m$. 
The precise forms of the last $m$ sums will be
specified below.

Let $\pbar_i=p''_i$ for $1\le i\le m$.
The choice of the $m$ index variables $\qbar_1,\dots,\qbar_m$
exactly corresponds to the choice of a $\bar{q}$ contributing
to (\ref{mss}) {\em together} with an
injection $\theta$, up to automorphisms of $q''$.
More precisely, as the index variables of the first $n$ sums determine
$$q=(q_1,\dots,q_n),$$ the unordered
partition $q'$ consists of $\qbar_1,\dots,\qbar_m$ and
$$d-\sum_{i=1}^n q_i-\sum_{j=1}^m 
\qbar_j$$ parts equal to $1$. Clearly, $q''$ consists
then of those $\qbar_j$ that are at least $2$. Finally, the $j$ such
that $\qbar_j\ge2$ determine $\im(\theta)$, and $\theta$ itself is
uniquely determined up to automorphisms of $q''$.

Let us analyze the factor $S[r''](q')$. Let 
$q^*$ denote 
the subpartition of $q'$ consisting of $\qbar_1,\dots,\qbar_m$.
Fix an injection 
$$\iota:\underline{\ell(r'')}\to \underline{\ell(q')}.$$
Suppose $c$ of the $m$ parts of $q^*$ lie in $\im(\iota)$.
The summand
$$ \prod_{j=1}^{\ell(r'')}\,\frac{1}{{q'_{\iota(j)}}^{r''_j-1}} \cdot
\prod_{i\notin\im(\iota)}\,\frac{1}{q'_i}$$
of $S[r''](q')$ corresponding to $\iota$ can, in first approximation, 
be accounted for by taking $m-c$ factors of type $(\beta)$ and
$c$ factors of type $(\gamma)$. 
The last $c$ factors correspond 
exactly to the $c$ parts of $q^*$ in $\im(\iota)$, and the constants
$\rbar$ in the factors are set to equal the parts $r''_j$.

Given $\iota$, we have defined an initial product of $n+m$ factors
of type $(\alpha)$, $(\beta)$, and $(\gamma)$.
However, the terms of the product do not properly account
for all the data of $\iota$.
The difficulty is that while $\iota$
specifies the subinjection of the $\ell(r'')-c$ parts of
$r''$ to the parts of $q'$ not in $q^*$, the
terms of the defined product do not.

We must alter the product in order to properly account for
$\iota$.
The partition $q'$ 
consists of $q^*$ and 
$$e=d-\sum_{i=1}^n q_i-\sum_{j=1}^m \qbar_j$$ 
parts equal to $1$.
The number of subinjections of the remaining $\ell(r'')-c$ parts of
$r''$ into the $e$ parts of $q'$ not in $q^*$ equals
$$e!/(e-\ell(r'')+c)!.$$ 
By giving such an subinjection in addition to the other data, the injection
$\iota$ is essentially determined.
To be precise, by the choice of
$c$ factors of type $(\gamma)$, the choice of $c$ corresponding parts
$r''_j$,
 and the choice of a subinjection of the remaining $\ell(r'')-c$ parts of
$r''$ into the $e$ parts of $q'$ not in $q^*$,
the injection $\iota$
is uniquely determined up to automorphisms of the $c$ chosen parts
$r''_j$.

Consider 
the $e!/(e-{\ell(r'')} +c)!$ 
subinjections into the $e$ parts of $q'$ not in $q^*$,
with 
$$e=d-\sum_{i=1}^n q_i-\sum_{j=1}^m \qbar_j.$$ 
A crucial observation is the following:
we can take account of the subinjections
by slight modifications of the factors of type $(\alpha)$, $(\beta)$,
and $(\gamma)$. Namely, after expanding $e!/(e-\ell(r'')+c)!$ we find
a polynomial $f$ of degree $\ell(r'')-c$ in $q_i$ and $\qbar_j$. Each
$q_i$ raises the exponent of $q_i$ in the $i$th factor of type $(\alpha)$
by $1$. Similarly, each $\qbar_j$ raises the exponent of $\qbar_j$
in the $j$th factor of type either $(\beta)$ or $(\gamma)$ by $1$.
After raising, a factor of type $(\beta)$ becomes 
a factor of type $(\beta')$.

Starting with an injection $\iota$, we have obtained, by expanding $f$,
a sum of products. Each product consists of factors (possibly raised)
of types $(\alpha)$, $(\beta)$, and $(\gamma)$.
The sum of products defines a
subsum of
$$
\sum_{\bar{q}}
\prod_{j=1}^n\,\Big(\binom{p_j-1}{q_j-1}(-1)^{q_j-1}q_j^{p_j-1-r_j}\Big)
\cdot \frac{T[q''](p'')}{|\Aut(q'')|}\cdot S[r''](q'). $$
For example, the factor $2^{v(\theta)}$ in $T[q''](p'')$
has been rewritten as a product
of the factors $2^{\delta(2,\pbar_j,\qbar_j)}$ occurring in the
(possibly raised) sums $(\beta)$ and $(\gamma)$.

We may now finally establish the required vanishing and non-vanishing.
Since $\bar{r}\le \bar{p}$, we have 
$$\ell(r'')\le\ell(p'')=m.$$ 
The initial product obtained from $\iota$
had $n$ factors of type $(\alpha)$, $m-c$ factors of type $(\beta)$,
and $c$ factors of type $(\gamma)$. Every factor of type $(\beta)$
vanishes. Every such factor must be raised for a non-vanishing
contribution. The maximum number of raisings is $\ell(r'')-c$. So for
non-vanishing, $\ell(r'')=m$ is required. Moreover, every factor of type $(\beta)$
becomes a factor of type $(\beta')$ with a minus sign. The product
of $n$ factors of type $(\alpha)$, $m-c$ factors of type $(\beta')$,
and $c$ factors of type $(\gamma)$, the unique
possibly non-vanishing contribution, appears with a coefficient
$(-1)^{m-c}(m-c)!$. 

Since $\ell(r'')=\ell(p'')$, the partition $r''$ must preceed $p''$
in the lexicographic order. So for every part of $p''$ equal
to $2$, there exists a part of $r''$ equal to $2$. To obtain
a non-vanishing contribution,  each such part
of $p''$ must come with a factor of type $(\beta')$. 
Considering the parts of $p''$ equal to 
$$3,4,5,\ldots$$ 
successively,
we find a non-vanishing contribution requires that each such part
come with a factor of type $(\gamma)$ and an equal part of $r''$.
Thus, $r''=p''$. Arguing as in the case $r''=p''=\emptyset$, we find
 a non-vanishing contribution requires $\bar{r}=\bar{p}$. We have proven
the required vanishing when $\bar{r}<\bar{p}$.

For the diagonal elements, assume
$\bar{r}=\bar{p}$. The only non-vanishing contributions arise by exactly
matching every part of $p'''$ with an equal part of $r''$
(via the element of $q'$ in the corresponding factor of type
$(\gamma)$), inserting a factor of type $(\beta')$ for
each remaining $2$ in $p''$, and mapping the $2$'s in $r''$ to the
$1$'s of $q'$ not in $q^*$.
The number of automorphisms
of the $c$ chosen parts $r''_j$ may now be identified with
$|\Aut(p''')|$. The factor $(m-c)!$ may now be identified
with the number of permutations of the $2$'s. Multiplied together,
these factors yield $|\Aut(p'')|$, cancelling the denominator
of the prefactor of the matrix element $\matc(\bar{p},\bar{p})$.
Finally, the factor $(-1)^{m-c}$ cancels the contribution
of the $m-c$ factors of type $(\beta')$, and the numerator
$\prod\,p_h\,\prod\,p'''_i$ of the leading factor cancels
the contribution of the factors of type $(\alpha)$ and $(\gamma)$.
We conclude $C(\bar{p},\bar{p})=1$ for all $\bar{p}$.
\epf

The lower bound $d-k$ for the length of
partitions $\bar{p}$ in $\Pi(d,n,k)$ does not play a crucial role.
In fact, 
Lemmas \ref{lmaa}
and \ref{lmbb} hold for the submatrices of $\mata(d,n,\infty)$ and 
$\matb(d,n,\infty)$
indexed by any subset $\Xi$ of $\Pi(d,n,\infty)$
for which
$$\Pi_{\bar{p}}=\{\,\bar{q}\in\Pi(d,n,\infty):\quad q\le p, 
\quad 1^{\ell(p'')-\ell(q'')}\,q''\le p''\,\}$$
is contained in $\Xi$ for all $\bar{p}\in\Xi$.
The proof is obtained simply by 
analyzing the vanishing of entries of
$\matb(d,n,\infty)$
forced by the definition of the matrix 
 (and the definition of $T$).

For Lemmas \ref{lmaa} and \ref{lmbb},
the set $\Xi$ should be closed under the operation of lowering
the parts of 
a partially ordered partition, necessarily making up for this by adding
parts equal to $1$ to the unordered partition. Clearly, 
the subsets $\Pi(d,n,k)$ satisfy this condition for all $k\ge0$.

\subsection{The nonsingularity of $\mata(\DD,\NN,k)$} 
The nonsingularity of $\mata(\DD,\NN,k)$ is immediate
when 
$$k\ge d-\sum_{i=1}^c|n_i|,$$ 
for then
the matrix is the Kronecker product of the matrices $$A(d_i,|n_i|,d_i-|n_i|).$$
For smaller $k$, certain rows and the corresponding columns are omitted
and the resulting matrix is not necessarily a Kronecker product.
However, for each $k$ the set $\Pi(\DD,\NN,k)$ is closed 
under the operation of lowering
the parts of some of the constituent 
partially ordered partitions and compensating
by adding parts equal to $1$ to the corresponding unordered partitions.
From the remark at the end of \ref{plmaa},  the method
of proof of Lemmas \ref{lmaa} and \ref{lmbb} still applies.

The nonsingularity of $\mata(\DD,\NN,k)$ is established
and the proof of Lemma \ref{mmm} is complete. 
\epf

\section{Consequences}

\subsection{Proof of Proposition \ref{jjjj}}
\label{vs}

\subsubsection{The moduli space $\overline{M}_{g,n}$} We prove the socle,
$$R^{3g-3+n}(\overline{M}_{g,n}),
$$ is $1$-dimensional. By [GrP2],
Proposition 11, the socle is generated by classes of the form
$$\xi_{B*}(\prod_{v\in V(B)} \theta_v),$$
where $B$ is a stable graph of genus $g$ with $n$ legs and 
$\theta_v\in R^*(\overline{M}_{g(v),n(v)})$ is a monomial in
the $\psi$ and $\kappa$ classes of $\overline{M}_{g(v),n(v)}$,
of degree equal to the dimension
$3g(v)-3+n(v)$. By Proposition \ref{jj}, we may assume 
$$3g(v)-3+n(v)\le g(v)-1+\delta_{0g(v)},$$
for all $v\in V(B)$. This implies $g(v)=0$ and $n(v)=3$.
So $R^{3g-3+n}(\overline{M}_{g,n})$ is generated by the point classes
of maximally degenerate stable $n$-pointed curves of genus $g$. All
such curves are degenerations of stable $n$-pointed irreducible rational
curves (with $g$ nodes). The closed stratum of such curves is
dominated by $\overline{M}_{0,2g+n}$. Hence all the point classes are
equivalent, which proves the socle claim for $\overline{M}_{g,n}$.

\subsubsection{The moduli space $M_{g,n}^c$} 
The maximally degenerate
stable $n$-pointed curves of genus $g$ of compact type consist of
$g$ elliptic tails attached to $g$ of the marked points of a
maximally degenerate stable $(g+n)$-pointed curve of genus $0$
(contracted to a point when $g+n=2$).
The classes of the different strata determined
by these maximally degerate curves are rationally equivalent in
$M_{g,n}^c$ via
the rational equivalence of points on $\overline{M}_{0,g+n}$.
We prove that these strata classes generate
the socle and that classes of higher degree vanish. 

The tautological ring
$R^*(M_{g,n}^c)$ is additively generated by classes of the form
$\xi_{B*}(\prod_{v\in V(B)} \theta_v)$, where $B$ is now a stable tree
of genus $g$ with $n$ legs. By Proposition \ref{jj}, we may assume
$$\deg \theta_v\le g(v)-1+\delta_{0g(v)}-\delta_{0n(v)}.$$
Since
$$ g(v)-1+\delta_{0g(v)}-\delta_{0n(v)}\le 2g(v)-3+n(v),$$
we find 
\begin{eqnarray*}
\deg \xi_{B*}(\prod_{v\in V(B)} \theta_v) &=& -1+\sum_{v\in V(B)}
(1+\deg \theta_v) \\
&\le& -1+\sum_{v\in V(B)}(g(v)+\delta_{0g(v)}-\delta_{0n(v)}) \\
&\le& -1+\sum_{v\in V(B)}(2g(v)-2+n(v)) = 2g-3+n,\\
\end{eqnarray*}
which proves the vanishing. The equality
$$ \deg \xi_{B*}(\prod_{v\in V(B)} \theta_v)=2g-3+n$$
implies that $(g(v),n(v))=(0,3)$ or $(1,1)$ for all $v$,
proving the socle claim for $M_{g,n}^c$.

\subsubsection{The moduli space $M_{g,n}^{rt}$ for $g\ge2$}
The  tautological ring is
additively generated by classes of the form
$\xi_{B*}(\prod_{v\in V(B)} \theta_v)$, where $B$ is now a stable tree
of genus $g$ with $n$ legs and a single vertex $w$ of genus $g$ (and all other
vertices of genus $0$). By Proposition \ref{jj}, we may assume
the class $\theta_v$ has degree
at most $g-1-\delta_{0n}$ on $w$ (and degree 0 on the genus 0 vertices). Then
\begin{eqnarray*}
\deg \xi_{B*}(\prod_{v\in V(B)} \theta_v) &=& -1+\sum_{v\in V(B)}
(1+\deg \theta_v) \\
&\le& g-1-\delta_{0n}-1+\sum_{v\in V(B)} 1 \\
&\le& g-1-\delta_{0n} + \sum_{v\neq w} (n(v)-2) \\
&=& g-1-\delta_{0n} + n-n(w) \le g-2+n,
\end{eqnarray*}
which proves the vanishing. The equality
$$ \deg \xi_{B*}(\prod_{v\in V(B)} \theta_v)=g-2+n$$
implies  $n(v)=3$ for $v\neq w$, that $n(w)=0$ or $1$, and 
$\theta_w$ is a top class on $M_g$ or $M_{g,1}$. The $1$-dimensionality
of the socle follows now from Looijenga's results 
[Lo] for $M_g$ and $M_{g,1}$,
the nonvanishing of $\kappa_{g-2}$ (see [FP]), and the rational equivalence
of points on $\overline{M}_{0,n+1}$.

\subsection{Gromov-Witten theory}
\label{gw}
\subsubsection{$\kappa$ descendent invariants}
Let $X$ be a nonsingular projective variety.
The $\kappa$ descendent Gromov-Witten invariants of $X$ are defined by:
\begin{equation}
\label{kejj}
\langle\ \tau_{e_1}(\gamma_1) 
\cdots \tau_{e_n}(\gamma_n) \prod_{j\geq 0} \kappa_j^{f_j}\
\rangle_{g,n,\beta}^X
= \int_{[\overline{M}_{g,n}(X,\beta)]^{vir}} 
\prod_{i=1}^n \psi_i^{e_i} \text{ev}_i^*(\gamma_i)
\prod_{j\geq 0} \kappa_j^{f_j},
\end{equation}
where $\psi_i$ is the cotangent line on the domain,
$\text{ev}_i^*(\gamma_i)$ is the pull-back of
$\gamma_i \in H^*(X,{\mathbb{Q}})$ via the $i^{th}$ evaluation map, and 
$\kappa_j$ is the Arbarello-Cornalba $\kappa$ class on the
       moduli space of maps.

\subsubsection{Proof of Proposition \ref{jjjjj}}
Since Proposition \ref{jjjjj} is 
well-known in genus 0 and 1, we will assume $g\geq 2$.
Let
$$\rho: \overline{M}_{g,n}(X,\beta) \rarr \overline{M}_{g,n}.$$
By the comparison results relating $\psi_i$ and $\kappa_j$
to the $\rho$-pull backs of the corresponding classes on the moduli
space of curves,
the invariants (\ref{kejj}) are equal to
the integrals 
\begin{equation}
\label{kejjj}
\int_{[\overline{M}_{g,n}(X,\beta)]^{vir}} 
\prod_{i=1}^n \rho^*(\psi_i^{e_i}) \text{ev}_i^*(\gamma_i)
\ \prod_{j\geq 0}\rho^*( \kappa_j^{f_j}),
\end{equation}
modulo corrections by invariants indexed by  lower data  $(g,n,\beta)$.
Here, $\psi_i$ and  $\kappa_j$ denote the cotangent line and
Arbarello-Cornalba $\kappa$ class
on the moduli space of curves. We will prove
 Proposition \ref{jjjjj} for the integrals (\ref{kejjj}). The
result for the invariants (\ref{kejj}) follows.

Let $\sum_{i=1}^n e_i + \sum_{j\geq 0} jf_j \geq g$.
Then, by Proposition \ref{jj}, the
class,
\begin{equation}
\label{pzx}
\prod_{i=1}^n \psi_i^{e_i} \prod_{j\geq 0} \kappa_j^{f_j},
\end{equation}
on $\overline{M}_{g,n}$
can be rewritten as a tautological boundary class.
The boundary strata are indexed by dual graphs.
Basic tautological classes are obtained on a 
given stratum by 
 products of $\psi$ and
$\kappa$ classes,
 $$\prod_{i=1}^{val(v)} \psi_i^{a_i} \prod_{j\geq 0} \kappa_j^{b_j},$$
at each vertex $v$ of the dual graph.
Every tautological boundary class is a linear combination of
these basic classes (see [GrP2]).
By Proposition  \ref{jj}, we may require the condition
\begin{equation}
\label{pzzx}
\sum_{i=1}^{val(v)} a_i + \sum_{j\geq 0} jb_j < g(v) + \delta_{0g(v)} 
\end{equation}
at each vertex $v$.

The proof of Proposition \ref{jjjjj} is completed by 
rewriting (\ref{pzx}) as a sum of basic classes satisfying the
vertex condition (\ref{pzzx}) and then applying the splitting
axiom to (\ref{kejjj}).
\epf

\subsection{The generation conjecture}
\label{lll}
We conjecture the ring $R^*(M_{g,n})$ is 
additively
generated 
by a restricted set of tautological classes.

\medskip
\noindent {\bf Conjecture 4.}
{\em $R^*(M_{g,n})$ is
additively generated 
by monomials
$$\prod_{i=1}^n \psi_i^{e_i} \prod_{j\geq 1} \kappa_j^{f_j}$$
satisfying
$$\sum_{i=1}^{n} e_i + \sum_{j\geq 1} (j+1){f_j} <g + \delta_{0g}.$$}

\medskip
We have verified Conjecture 4 for $g+3n\leq 21$, assuming for $n>0$
the Gorenstein conjecture for ${\mathcal C}_g^n$, the fiber product of the
universal curve over $M_g$.

Conjecture 3 is a consequence of the following three
statements:
 Proposition \ref{jjjjj},
Conjecture 4, and the (unproven)
exactness of the sequence
$$R^*(\partial \overline{M}_{g,n}) \rarr
R^*(\overline{M}_{g,n}) \rarr R^*(M_{g,n}) \rarr 0.$$
Conjecture 3 is derived from these statements
 by the standard method of expressing
$\kappa$ monomials in terms of push-forwards of 
cotangent line classes on additionally pointed spaces.
On $\overline{M}_{g,n}$, the class,
$$\prod_{j\geq 1} \kappa_j^{f_j},$$
is expressed as
a sum of classes,
$$\pi_{r*} (\prod_{i=1}^r \psi_{n+i}^{h_i} ),$$
where
$$\pi_r: \overline{M}_{g,n+r} \rarr \overline{M}_{g,n}$$
is the forgetting map and
$$r\leq \sum_{j\geq 1} f_j, \qquad
 \sum_{i=1}^r h_i \leq \sum_{j\geq 1} (j+1)f_j,$$
see [AC].
We leave the details to the reader.

%\vspace{+10 pt}
%\noindent Institutionen f\"or Matematik \\
%\noindent Kungl Tekniska H\"ogskolan \\
%\noindent 100 44 Stockholm, Sweden \\
%\noindent faber@math.kth.se

%\vspace{+10 pt}
%\noindent Department of Mathematics\\
%\noindent Princeton University\\
%\noindent Princeton, NJ 08544\\
%\noindent rahulp@math.princeton.edu
%\end{document}

\vspace{+10 pt}
\noindent Institutionen f\"or Matematik \hfill Department of Mathematics\\
\noindent Kungl Tekniska H\"ogskolan \hfill Princeton University\\
\noindent 100 44 Stockholm, Sweden \hfill Princeton, NJ 08544\\
\noindent faber@math.kth.se \hfill rahulp@math.princeton.edu

\begin{thebibliography}{[COGP]}

\bibitem[AC] {ac} E.
Arbarello and M. Cornalba, {\em 
Combinatorial and algebro-geometric cohomology classes on the moduli
   spaces of curves.} J. Algebraic Geom. {\bf 5} (1996), no. 4, 705--749. 

\bibitem[BP]{bp} P. Belorousski and R. Pandharipande, {\em 
A descendent relation in genus 2}.
Ann. Scuola Norm. Sup. Pisa Cl. Sci. (4) {\bf 29} (2000), no. 1, 171--191.

\bibitem[FP]{fp} C. Faber and R. Pandharipande, {\em
Logarithmic series and Hodge integrals in the tautological ring}.
With an appendix by Don Zagier.
Michigan Math. J. {\bf 48} (2000), 215--252.

\bibitem[Ge1]{ge1} E. Getzler, {\em
Intersection theory on $\overline{M}_{1,4}$ and elliptic Gromov-Witten 
invariants}.  J. Amer. Math. Soc. {\bf 10}  (1997),  no. 4, 973--998. 

\bibitem[Ge2]{ge2} E. Getzler, {\em
Topological recursion relations in genus $2$}.  Integrable systems and 
algebraic geometry (Kobe/Kyoto, 1997),  73--106, World Sci. Publishing, 
River Edge, NJ, 1998.

\bibitem[GrP1]{grp1} T. Graber and R. Pandharipande, {\em
Localization of virtual classes}.  
Invent. Math.  {\bf 135}  (1999),  no. 2, 487--518.

\bibitem[GrP2]{grp2} T. Graber and R. Pandharipande, {\em
Constructions of nontautological classes on moduli spaces of curves}.
math.AG/0104057.

\bibitem[GrV1]{grv1} T. Graber and R. Vakil, {\em
On the tautological ring of $\overline{M}_{g,n}$}.  
Turkish J. Math. {\bf 25}  (2001),  no. 1, 237--243.

\bibitem[GrV2]{grv2} T. Graber and R. Vakil, 
in preparation.

\bibitem[HM]{hm} J. Harris and D. Mumford, {\em
On the Kodaira dimension of the moduli space of curves}. 
With an appendix by William Fulton.  
Invent. Math.  {\bf 67}  (1982), no. 1, 23--88.

\bibitem[I]{i} E.-N. Ionel, {\em
Topological recursive relations in $H^{2g}(M_{g,n})$}.  
Invent. Math.  {\bf 148}  (2002),  no. 3, 627--658.

\bibitem[K]{k} S. Keel, {\em
Intersection theory of moduli space of 
stable $n$-pointed curves of genus zero}.  
Trans. Amer. Math. Soc.  {\bf 330}  (1992),  no. 2, 545--574.

\bibitem[Li]{li} J. Li, {\em
Stable morphisms to singular schemes and relative stable morphisms}.  
J. Differential Geom. {\bf 57}  (2001),  no. 3, 509--578.

\bibitem[LR]{lr} A.-M. Li and Y. Ruan, {\em
Symplectic surgery and Gromov-Witten invariants of Calabi-Yau 3-folds}.  
Invent. Math. {\bf 145}  (2001),  no. 1, 151--218.

\bibitem[Lo]{lo} E. Looijenga, {\em
On the tautological ring of $M_g$}.  
Invent. Math. {\bf 121}  (1995),  no. 2, 411--419.

\bibitem[P1]{p1} R. Pandharipande, {\em
A geometric construction of Getzler's elliptic relation}.  
Math. Ann. {\bf 313}  (1999),  no. 4, 715--729.

\bibitem[P2]{p2} R. Pandharipande, {\em
Three questions in Gromov-Witten theory}.
Proceedings of the ICM 2002 Beijing, Vol. II, 503--512.  


\end{thebibliography}
\end{document}